\pdfoutput=1
\documentclass{article}
\usepackage[T1]{fontenc}
\usepackage[utf8]{inputenc}

\usepackage[final]{microtype}

\usepackage{mathtools}
\mathtoolsset{showonlyrefs}

\usepackage{amsthm}
\theoremstyle{plain}
\newtheorem{theo}{Theorem}[section]
\newtheorem{theorem}{Theorem}[section]

\newtheorem{coro}{Corollary}[section]
\newtheorem{corollary}{Corollary}[section]

\newtheorem{lemma}{Lemma}[section]
\theoremstyle{definition}
\newtheorem{defi}{Definition}[section]
\theoremstyle{remark}

\usepackage{amssymb}
\usepackage{latexsym}
\usepackage{epsfig}
\usepackage{graphics}
\usepackage{array}
\usepackage{euscript}
\usepackage{dsfont}
\usepackage{esint}
\usepackage{slashed} 
\usepackage{faktor} 

\usepackage{tikz-cd} 
\tikzset{
	pf/.style={commutative diagrams/.cd, every arrow, every label},
	surj/.style=commutative diagrams/two heads,
	inj/.style=commutative diagrams/hook,
	gl/.style=commutative diagrams/equal,
	mat/.style={matrix of math nodes, commutative diagrams/.cd, every cell},
	dr/.style={matrix of math nodes, commutative diagrams/.cd, every cell, column sep=small},
	seq/.style={matrix of math nodes, commutative diagrams/.cd, every cell, column sep=small}
	}
\newenvironment{diag*}{\[\begin{tikzpicture}[commutative diagrams/.cd, every diagram, baseline=(current bounding box.center)]}{\end{tikzpicture}\]\ignorespacesafterend}
\newenvironment{diag}{\begin{equation}\begin{tikzpicture}[commutative diagrams/.cd, every diagram, baseline=(current bounding box.center)]}{\end{tikzpicture}\end{equation}\ignorespacesafterend}

\usepackage[final]{hyperref}

\newcommand{\Z}{\mathbb{Z}}
\newcommand{\R}{\mathbb{R}}
\newcommand{\C}{\mathbb{C}}

\DeclareMathOperator{\ima}{i}
\DeclareMathOperator{\PSL}{PSL}
\DeclareMathOperator{\Diff}{Diff}

\newcommand{\CliffordAlgebraM}{\operatorname{Cl}_{0,2}}
\newcommand{\CliffordAlgebraP}{\operatorname{Cl}_{2,0}}
\newcommand{\CliffordBundleM}{{\operatorname{Cl}\left(\Sigma, -\gamma\right)}}
\newcommand{\CliffordBundleP}{{\operatorname{Cl}\left(\Sigma, \gamma\right)}}

\newcommand{\SpinGroup}{\operatorname{Spin}(2)}
\newcommand{\SpinorBundleM}{S}
\newcommand{\SpinorBundleP}{\mathcal{S}}

\newcommand{\Smooth}[1]{{|#1|}}
\newcommand{\SpinStructure}{P_{\SpinGroup}}
\newcommand{\Rm}{R}

\newcommand{\ACRS}{J_\Sigma}
\newcommand{\ACSBP}{J_\SpinorBundleP}

\newcommand{\dirac}{\slashed{\partial}}
\newcommand{\Dirac}{\slashed{D}}

\newcommand{\dd}{\mathop{}\!\mathrm{d}}
\newcommand{\dv}{\dd{vol}}
\newcommand{\dx}{\dd{x}}
\newcommand{\dz}{\dd{z}}
\newcommand{\dy}{\dd{y}}

\newcommand{\Tp}{\mathcal{T}_p}
\newcommand{\ModuliSpace}{\mathcal{M}_p}

\newcommand{\ADH}{A_{DH}} 
\newcommand{\ADHG}{\mathbb{A}} 
\newcommand{\ASRS}{A} 

\newcommand{\bel}[1]{\begin{equation}\label{#1}}
\newcommand{\be}{\begin{equation}}
\newcommand{\qe}{\end{equation}}

\newcommand{\rf}[1]{\eqref{#1}}

\newcommand{\p}{\partial}

\newcommand{\na}{\nabla}
\newcommand{\al}{\alpha}

\newcommand{\cD}{\mathcal{D}}
\newcommand{\cO}{\mathcal{O}}

\DeclareMathOperator{\trace}{trace}

\DeclareMathOperator{\diverg}{div}

\title{From harmonic maps to the nonlinear supersymmetric sigma model of quantum field theory.\\
At the interface of theoretical physics, Riemannian geometry and nonlinear analysis}
\author{Jürgen Jost \and Enno Keßler \and Jürgen Tolksdorf \and Ruijun Wu \and Miaomiao Zhu}

\date{}

\begin{document}
\maketitle
\begin{center}
	\textbf{\large  Dedicated to the memory of Eberhard Zeidler}
\end{center}

\begin{abstract}
  Harmonic maps from Riemann surfaces arise from a conformally invariant variational problem. Therefore, on one hand, they are intimately connected with moduli spaces of Riemann surfaces, and on the other hand, because the conformal group is noncompact, constitute a prototype for the formation of singularities, the so-called bubbles, in geometric analysis. In theoretical physics, they arise from the nonlinear $\sigma$-model of quantum field theory. That model possesses a supersymmetric extension, coupling a harmonic map like field with a nonlinear spinor field. In the physical model, that spinor field is anticommuting. In this contribution, we analyze both a mathematical version with a commuting spinor field and the original supersymmetric version. Moreover, this model gives rise to a further field, a gravitino, that can be seen as the supersymmetric partner of a Riemann surface metric. Altogether, this leads to a beautiful combination of concepts from quantum field theory, structures from Riemannian geometry and Riemann surface theory, and methods of nonlinear geometric analysis. 
\end{abstract}

\section{Introduction}
At his premature death, Eberhard Zeidler left his magnum opum unfinished, a comprehensive mathematical treatment of quantum field theory.
He had planned six volumes, of which three, with more than 3000 pages,~\cite{Zeidler06,Zeidler09,Zeidler11}, could be realized, and it is planned that one of us (J.T.) can complete the fourth volume on the basis of the notes that he had left behind.
His vision consisted into a profound interpenetration of the two fields, leading to deeper insight and new methods in both of them.

Indeed, the action functionals of (quantum) field theory usually contain fundamental geometric structures that are incorporated in their symmetries and invariances.
In turn, analyzing their Euler--Lagrange equations---or physically speaking, field equations---from a rigorous mathematical perspective leads to challenging problems in nonlinear partial differential equations and the calculus of variations.
And solving those problems can provide us with insights and methods that can find profound novel applications in geometry, as for instance the Yang-Mills and Seiberg-Witten equations have shown.
And superstring theory, while controversial among physicists, has triggered important advances in various areas of mathematics.

In this spirit, in this contribution, we discuss some of the basic models of quantum field theory, the nonlinear $\sigma$-model.
In its most basic version, mathematically, it involves harmonic mappings from a Riemann surface to a sphere or some other model space.
Since the target is not Euclidean, but rather a symmetric space, that is, a particular Riemannian manifold, the Euler--Lagrange equations of the model are generically nonlinear.
Their mathematical treatment uncovered the phenomenon of bubbling of solutions of conformally invariant variational problems, ultimately leading to the development of the mathematical theory of pseudoholomorphic curves and Gromov-Witten invariants. Also, there are many analogies to the four-dimensional Yang-Mills equations.

From the perspective of physics, the ultimate version of the model is still richer, as it includes supersymmetry.
That is, the original bosonic fields are supplemented by fermionic partners.
This brings us into the realm of supergeometry and, in particular, the theory of super Riemann surfaces.

Here, on one hand, we confront that challenge and develop the mathematical theory of the supersymmetric nonlinear $\sigma$-model. As we shall see, the geometric perspective, that is, the moduli space of super Riemann surfaces, and the physical or variational one, the symmetries and invariances of the action functional, are two sides of the same coin, and either of them deeply reflects the other. More precisely, the Noether currents associated to the invariances of the functional yield the cotangent vectors of the moduli space of super Riemann surfaces.

We have, however, discovered more. There is also a version of the theory that does not need supergeometry. This converts the formal expansion that supergeometry works with into nonlinear elliptic systems of partial differential equations that extend the harmonic map system. These systems on one hand provide us with a test bed for the development of new analytical methods that should eventually find applications in the wider realm of geometric analysis, and on the other hand uncover some novel geometric structure that seems to be worth in exploring further, in order to create new types of invariants of Riemannian manifolds.

We gratefully recall the inspiration that Eberhard Zeidler has given us over the years, through his works, his generous and warm personality, and his scientific vision, and we hope that it is not inappropriate to dedicate the present work to his memory.

\section{The geometric and analytic context: Harmonic mappings between Riemannian manifolds}\label{harmonic}
Harmonic mappings from Riemann surfaces into Riemannian manifolds can be considered as special cases of harmonic mappings between Riemannian manifolds.
While, as we shall see and explore, when the domain is a Riemann surface, harmonic mappings enjoy special properties, on which the rest of this article will essentially build, it is nevertheless useful to first understand the wider context.
Thus, let $M$ be a Riemannian manifold of dimension $m$.
In the sequel, we shall use basic concepts of Riemannian geometry without further explanation; for example, $TM$ is the tangent bundle of $M$, and $T_xM$ is the tangent space over the point $x\in M$.
The reader may consult~\cite{Jost17}.

We shall usually implicitly assume that all manifolds are compact, without boundary. While boundary value problems are also analytically challenging and geometrically interesting, for this article, we have decided not to go into that topic.

We need some notation. In local coordinates, the metric tensor of $M$ is written as
\begin{equation}%
\label{ha13}
	{\left(\gamma_{\alpha \beta}\right)}_{\alpha, \beta = 1, \dotsc, m}.
\end{equation}
More precisely, the Riemannian metric is the symmetric 2-tensor
\begin{equation}
\label{ha14}
	\gamma=\gamma_{\alpha \beta}\dx^\alpha \otimes \dx^\beta,
\end{equation}
but for most purposes, we shall just work with its coefficient matrix \rf{ha13}.
From~\rf{ha14}, the reader can already see that we use the standard Einstein summation convention.
We shall also use the following notations
\begin{alignat}{2}\label{ha15}
	& {\left(\gamma^{\alpha\beta}\right)}_{\alpha,\beta=1,\dotsc , m} = {\left(\gamma_{\alpha \beta}\right)}^{-1}_{\alpha , \beta} &
		\quad & \text{(inverse metric tensor)}\\
	& |\gamma| \coloneqq \det (\gamma_{\alpha\beta}) & &\\
\label{ha16}	& \Gamma^\alpha_{\beta \eta} \coloneqq \frac{1}{2}\gamma^{\alpha \delta}
		(\gamma_{\beta \delta , \eta} + \gamma_{\eta \delta , \beta} - \gamma_{\beta \eta, \delta}) & &
		\text{(Christoffel symbols of $M$)}.
\end{alignat}
The volume form of $M$, $\sqrt{|\gamma|}\dx^1\cdots \dx^m$ in local coordinates, will often be abbreviated as $\dv_\gamma$. The Laplace-Beltrami operator of $M$ is
\begin{equation}\label{eq020112}
	\Delta_M  \coloneqq  \frac{1}{\sqrt{|\gamma|}}\frac{\p}{\p x^\alpha}\left( \sqrt{|\gamma|} \gamma^{\alpha \beta} \frac{\p}{\p x^\beta}\right).
\end{equation}
We employ here a sign convention different from that in \cite{Jost17}, where $\Delta_M$ has a minus sign, in order to make it a positive operator.

\subsection{Harmonic functions}\label{func}
For a function $\phi:M\to \R$, we  have the Dirichlet integral\index{Dirichlet integral}
\bel{ha20}
E(\phi) \coloneqq \frac{1}{2} \int_M \Vert d\phi\Vert^2 \dv_\gamma \coloneqq \frac{1}{2} \int_M \gamma^{\alpha\beta}(x) \frac{\partial \phi}{\partial x^\alpha} \frac{\partial \phi}{\partial x^\beta}
\sqrt{|\gamma|} \dx^1 \ldots  \dx^m.
\qe
Here, we have used an invariant and a coordinate notation, that is,
\bel{ha20a}
\Vert d\phi\Vert^2=\gamma^{\alpha\beta}(x) \frac{\partial\phi}{\partial x^\alpha} \frac{\partial \phi}{\partial x^\beta}.
\qe
The critical points of $E$ are the harmonic functions of $M$, that is, the solutions of the corresponding Euler--Lagrange equation
\bel{ha21}
\Delta_M \phi =0.
\qe
For the relevant results in the calculus of variations, we refer to \cite{Jost/Li98}, and for the PDE aspects, to \cite{Jost13}. We need to recall a few facts here, however.
The Dirichlet integral \rf{ha20} is naturally defined on the Sobolev space $W^{1,2}(M)$ of functions that possess square integrable weak derivatives, but need not be any smoother than that. In fact, they need not even be continuous, let alone differentiable. Therefore, for such a Sobolev function, a priori \rf{ha21} need not make sense. Thus, from minimizing $E$, one first obtains a so-called weak solution, that is, a $\phi$ solving
\bel{ha22}
\frac{1}{2} \int_M \gamma^{\alpha\beta}(x) \frac{\partial\phi}{\partial x^\alpha} \frac{\partial \eta}{\partial x^\beta}
\sqrt{|\gamma|} dx^1 \ldots  dx^m =0
\qe
for all test functions $\eta \in W^{1,2}(M)$.
One then shows, and this is the prototype of elliptic regularity theory, that any solution of \rf{ha22} is in fact smooth and therefore solves \rf{ha21}, as one finds by integrating \rf{ha22} by parts to move the first derivative from $\eta$, which is allowed when $\phi$ is smooth.
We refer to~\cite{Jost13} for the systematic theory. Here, smoothness of a weak solution is easy because \rf{ha21} is a \emph{linear} elliptic partial differential equation.
Of course, this assumes that the metric $\gamma$ itself is smooth and bounded, and also that $\sqrt{|\gamma|}$ is bounded from below, so that the inverse metric tensor $\gamma^{\alpha \beta}$ is likewise smooth and bounded.
More precisely, one needs the uniform ellipticity condition
\begin{equation}
\label{ha23}
	\lambda \Vert \xi \Vert^2 \le \sqrt{|\gamma|} \gamma^{\alpha \beta} \xi_\alpha \xi_\beta \le \Lambda \Vert \xi \Vert^2
	\text{ for all } \xi=(\xi_1,\dots \xi_m)\text{ and some }0<\lambda,\Lambda <\infty.
\end{equation}
When this holds, and the quantities $\sqrt{|\gamma|} \gamma^{\alpha \beta}$ are measurable, then by the De Giorgi-Nash theorem, a weak solution of \rf{ha22} is H\"older continuous.
And if these coefficients possess more regularity, then so does the solution $\phi$.
In fact, by Schauder's theorem, if they are of class $C^{k,\alpha}$ for some $k\in \mathbb{N}, 0< \alpha <1$, then~$\phi \in C^{k+2,\alpha}$.
Proofs of these results can be found in~\cite{Jost13}.
Anticipating the key discussion of later sections, however, we not only want to treat the function~$\phi$ as a variable, but also the metric~$\gamma$.
Even though later, for two-dimensional domains, which are those of main interest in this article, we shall also derive equations for variations of the Dirichlet integral w.r.t.\ the metric~$\gamma$, these will not be differential equations that will yield any regularity for the metric.
Rather, there is some gauge invariance principle.
That principle says that, in dimension two, again, we can always choose local coordinates so that the metric becomes the Euclidean metric, that is, constant, and therefore as smooth as one wants.
And then $\phi$ is also smooth to any order, in fact even analytic.
However, when we still move further, to the supersymmetric sigma model, we get further variables, a spinor field that is the partner of the function~$\phi$, and a gravitino field that is the partner of the metric $\gamma$.
And then, subtle regularity issues will arise.
But before that, we need to generalize things in another direction.

\subsection{Harmonic mappings}
We shall now move to \emph{nonlinear} elliptic systems where the regularity for $\phi$  is no longer so easy, and in fact, not even true in general.

We move  into this nonlinear realm  by replacing the target $\R$ by another Riemannian manifold \((N,g)\) of dimension $n$, with metric tensor $g=g_{ij}\dx^i\otimes\dx^j$.
Analogously to \rf{ha15}, \rf{ha16}, we define \(g^{ij}\) to be the coefficients of the inverse metric and the Christoffel symbols \(\Gamma^i_{jk}\).

For a map  $\phi\colon M \to N$ in the appropriate Sobolev space $W^{1,2}(M,N)$ (see~\cite{Jost17} for the precise definition), we define its energy\footnote{One should rather call that expression \emph{action} instead of energy, in view of the connections with quantum field theory developed below, but apparently, the terminology was introduced by people who did not see the relation with physical action principles.}
of  $\phi\colon M \to N$ is
\begin{align}\label{eq080106}
	E(\phi)&\coloneqq \frac{1}{2}\int_M \Vert d\phi\Vert^2 \dv_\gamma \\
	&= \frac{1}{2}\int_M \gamma^{\alpha \beta}(x) g_{ij}(\phi(x))\frac{\partial \phi^i(x)}{\partial x^\alpha} \frac{\partial \phi^j(x)}{ \partial x^\beta}\sqrt{|\gamma|} \dx^1 \dots \dx^m
\end{align}
where in the last expression, we use local coordinates $(x^1,\ldots , x^m)$ on $M$ and $(\phi^1, \ldots , \phi^n)$ on $N$.
Analogously to \rf{ha20a},
\bel{ha27}
\Vert  \dd\phi\Vert^2= \gamma^{\alpha \beta}(x) g_{ij}(\phi(x))
		\frac{\partial \phi^i(x)}{\partial x^\alpha} \frac{\partial \phi^j(x)}{\partial x^\beta}.
\qe
The energy integral thus is defined in terms of the Riemannian metrics of $M$ and $N$. When $N=\R$, this reduces to \rf{ha20}.

By a somewhat tedious, but otherwise straightforward, computation, we can generalize \rf{ha21}.
\begin{lemma}\label{l080101}
The Euler--Lagrange equations\index{Euler--Lagrange
equation} for $E$ are
\begin{equation}\label{eq080107}
	\frac{1}{\sqrt{|\gamma|}} \frac{\partial}{\partial x^\alpha}
		(\sqrt{|\gamma|} \gamma^{\alpha \beta} \frac{\partial}{\partial x^\beta} \phi^i) +
		\gamma^{\alpha \beta} (x) \Gamma^i_{jk} (\phi(x)) \frac{\partial}{\partial x^\alpha} \phi^j
		\frac{\partial}{\partial x^\beta} \phi^k = 0
\end{equation}
for $i=1,\dots ,n$.
\end{lemma}
We usually write this system of equations as
\bel{ha30}
\tau(\phi) =0,
\qe
and call $\tau(\phi)$ the tension field of $\phi$.
\begin{defi}\label{d080102}
	Solutions of~\eqref{eq080107} (or, equivalently, of~\eqref{ha30}) are called \emph{harmonic maps}.\index{harmonic map}
\end{defi}
The Euler--Lagrange equation~\rf{eq080107} constitutes a semilinear elliptic system.
It is in general nonlinear, because of the geometry of~$N$ that enters through the Christoffel symbols~$\Gamma^i_{jk} (\phi(x))$ that are contracted with first derivatives of~$\phi$, resulting in an expression that is quadratic in those derivatives.
Therefore, weak solutions, defined analogously to~\rf{ha22}, in general need not be regular.
Since the regularity theory becomes much better when the domain is two-dimensional, our exclusive case of interest below, we don't enter here into that general theory, however.

As we have seen, harmonic maps are generalization of harmonic functions.
Both the Laplace-Beltrami operator $\Delta_M$ and its generalization, the tension field~$\tau$, involve the Levi-Civita connection of~$M$, the latter also that of~$N$.

The concept of harmonic maps was first proposed by S.~Bochner~\cite{Bochner40} in 1940.
See also~\cite{Fuller54}.
Subsequently, it was developed very extensively, and in fact too extensively to be reviewed here.
We refer the reader to~\cite{Jost08c} for a survey.

In particular, one may apply variational methods to the action functional~$E$. As an alternative, one may use a parabolic method, that is, study the flow
\bel{ha40}
\frac{\partial \phi}{\partial t}=\tau(\phi)
\qe
where $\phi$ now is a map from $M\times [0,\infty)$ to~$N$. Of course, one needs to show that a solution exists for all $t\in [0,\infty)$ and converges to a solution of \rf{eq080107} as $t\to \infty$. This is not always true, but holds when the target space $N$ has nonpositive curvature. Like the regularity theory for \rf{ha30}, the subject has been too extensively studied to be properly reviewed here.

\section{Harmonic maps from Riemann surfaces and conformal invariance}%
\label{harmconf}
In this section and the following, we collect some material that is discussed in more detail in~\cite{Jost91,Jost97c,Jost06,Jost08c,Jost09a,Jost17,Jost/Yau10}, where also further references can be found.
The energy integral, and therefore harmonic maps, have some special property when the dimension of the domain is two, because of conformal invariance, as we are now going to explain.
We denote the domain by~$\Sigma$ and let it be an oriented two-dimensional Riemannian manifold with metric $\gamma_{\alpha \beta}$.
In the sequel, we shall use different versions of the  uniformization theorem, and we want to describe them now.
First, there is the local uniformization theorem, saying that locally, any metric on a surface can be transformed into the Euclidean one by a suitable coordinate transformation.
This was first proved by Gauss for smooth metrics.
It was an important step towards the modern regularity theory of solutions of elliptic partial differential equations when Lichtenstein~\cite{lichtenstein16} proved such a result for metrics of class $C^{1,\alpha}$.
Lavrent'ev~\cite{Lavrentev35} then showed such a result for continuous metrics, and Morrey~\cite{Morrey38} obtained it for bounded ones; of course, all metrics have to be uniformly positive definite.
See also the presentation in~\cite{Jost84}.
A global uniformization theorem was first stated by Riemann, and it is known as the Riemann mapping theorem.
The final version was proved by Koebe who showed that any non-empty simply connected open subset $U$ of the complex  plane $\C$ other than $\C$ itself is conformally equivalent to the open unit disk.
That means that there exists a bijective conformal map from $U$ onto the unit disk.
Finally, there is Poincar\'e's uniformization theorem which says that the universal cover of a compact Riemann surface is conformally equivalent to either the unit sphere (for genus 0), the complex plane (for genus 1), or the hyperbolic plane (for genus $p \ge 2$).

By the uniformization theorem of Gauss, $\Sigma$ then obtains the structure of a Riemann surface, that is, we can locally find holomorphic coordinates $z=x+\ima y$ on $\Sigma$ such that the metric tensor \rf{ha13} becomes
\begin{equation}%
\label{rs1}
	\lambda^2 (z)\, \dz \otimes \dd\overline{z}
\end{equation}
where $\lambda(z)$ is a real-valued positive function.
Here, $\dz=\dx+ \ima\dy$, $\dd\overline{z}=\dx-\ima\dy$. Conversely, given a compact Riemann surface $\Sigma$, by the Poincar\'e uniformization theorem, we can  find a conformal metric as in \rf{rs1} with constant curvature; we shall be mainly interested in the case of genus $\ge 2$, the hyperbolic case, where this curvature is negative, $-1$.
This metric is unique (up to diffeomorphism, see below).
As such a surface $\Sigma$ is of the form $H/\Gamma$, where~$H\coloneqq\{x+\ima y\in\C :y>0\}$ is the upper half plane and $\Gamma$ is a discrete subgroup of the isometry group~$\PSL(2,\R)$ of $H$.
Here the group \(\PSL(2,\R)\) acts on \(H\) by Möbius transformations.
The metric on $S$ descends from the hyperbolic metric $\frac{1}{y^2}\dz\otimes\dd\overline{z}$ on $H$.

As before, let~$N$ be a Riemannian manifold with metric tensor $(g_{ij})$.
The energy \rf{eq080106} of a map $\phi\colon \Sigma \to N$ then becomes
\begin{equation}\label{eq080211}
\begin{split}
	E(\phi) & = \frac{1}{2}\int_\Sigma \frac{4}{\lambda^2(z)} g_{ij} \frac{\partial \phi^i}{\partial z}
		\frac{\partial \phi^j}{\partial \overline{z}} \frac{i}{2}\lambda^2(z)\dz \land \dd\overline{z}
		\\
	& = \phantom{\frac{1}{2}} \int_\Sigma g_{ij} \frac{\partial \phi^i}{\partial z}\frac{\partial \phi^j}{\partial \overline{z}}
		\ima\dz \land \dd\overline{z} .
\end{split}
\end{equation}
The last equation implies that the energy\index{energy}
of a map from a Riemann surface $\Sigma$ into a Riemannian manifold is \emph{conformally invariant} in the sense
that it does not depend on the choice of a metric on~$\Sigma$, but only on the Riemann surface\index{Riemann surface} structure.
This only holds when the dimension of the domain is 2, that is, when we look at maps from a surface.
The reason is that the volume form behaves like an $\frac{n}{2}$-th power of the metric, and therefore, only for~$n=2$, this can cancel the effect of the inverse metric tensor in the energy integral. Nevertheless, as we shall explore below, it turns out to be useful to keep the domain metric $\lambda^2$ in the picture, even though, by the observation we have just made, this introduces some redundancy.

Also, if~$k\colon\Sigma \to S$ is a bijective holomorphic or antiholomorphic map between Riemann surfaces then for any $\phi\colon S \to N$ (of class $C^1$)
\bel{rs2}
E(\phi \circ k) = E(\phi),
\qe
and if $\phi$ is harmonic, then so is $\phi\circ k$.

The harmonic map equation \rf{l080101} is then also independent of the choice of conformal metric on the domain.
This means that the map $\phi\colon\Sigma\to N$ of class $C^2$ is harmonic\index{harmonic} iff
\begin{equation}\label{eq080210}
	\frac{\partial^2 \phi^i}{\partial z \partial \overline{z}} + \Gamma^i_{jk}(\phi(z)) \frac{\partial\phi^j}{\partial z}
		\frac{\partial \phi^k}{\partial \overline{z}} = 0 \quad
		\text{for}\ i = 1, \ldots, \dim N.
\end{equation}
When the target is also a Riemann surface $S$ with metric
\bel{eq080210a}
\rho^2(\phi)\dd\phi \otimes \dd\overline{\phi},
\qe
then the harmonic map equation \rf{eq080210} becomes
\bel{eq080210b}
\phi_{z\overline{z}}+ \frac{2\rho_\phi}{\rho}\phi_z \phi_{\overline{z}}=0.
\qe
Thus, holomorphic or antiholomorphic maps between Riemann surfaces are harmonic as they obviously satisfy \eqref{eq080210b}.
In the sequel, we shall write~$\pm$~holomorphic to mean holomorphic or antiholomorphic.
The converse does not necessarily hold.
Harmonic maps exist in any homotopy class of maps between compact Riemann surfaces (unless the target is a sphere, that is, of genus 0), but~$\pm$~holomorphic maps $h\colon \Sigma\to S$ have to satisfy the Riemann-Hurwitz formula
\begin{equation}%
\label{RH1}
	2-2p=|m|(2-2p)-v_h,
\end{equation}
where $p$, \(q\) are the genera of $\Sigma$, \(S\), $m$ is the degree of $h$, and $v_h\ge 0$ is its total ramification index.
In particular, as a necessary topological condition for a harmonic map to be $\pm$ holomorphic, we have the inequality
\begin{equation}%
\label{RH2}
	\chi(\Sigma)\le |m|\chi(S),
\end{equation}
where $\chi$ denoteos the Euler characteristic.
In fact, even if~\rf{RH2} holds, in general there will exist no holomorphic map between $\Sigma$ and $S$, as this requires that the conformal structure of \(\Sigma\) be a conformal cover of the one of $S$.
Thus, in general, harmonic maps between Riemann surfaces \(\Sigma\) and \(S\) cannot be~$\pm$~holomorphic.
Also, as we see from~\rf{eq080210}, whether a map is harmonic depends on the metric of the target $S$ whereas the property of being $\pm$ holomorphic only depends on the conformal structure of $S$, but not on its metric.
Nevertheless, harmonic maps into Riemann surfaces enjoy special properties, as demonstrated by the following result of Schoen-Yau~\cite{Schoen/Yau78} and Sampson~\cite{Sampson78}
\begin{theorem}%
\label{diffeo}
	Let $\phi\colon \Sigma \to S$ be a harmonic map with $p=q$, $|\deg \phi|=1$.
	If the curvature \(K_2\) of $S$ satisfies $K_2\le 0$, then $\phi$ is a diffeomorphism.
\end{theorem}
(In fact, harmonic diffeomorphism exist even without the curvature restriction, see~\cite{Jost81,Jost/Schoen82,Coron/Helein89}.)
A key ingredient in the proof are the following identities.
We put, for a harmonic map $\phi\colon \Sigma\to S$ between Riemann surfaces with curvatures $K_1,K_2$,
\begin{align}%
\label{rs110}
	H &\coloneqq|\partial \phi|^2\coloneqq\frac{\rho^2}{\lambda^2}\phi_z \overline{\phi}_{\overline{z}} &
	L &\coloneqq|\overline{\partial} \phi|^2\coloneqq\frac{\rho^2}{\lambda^2}\phi_{\overline{z}}\overline{\phi}_z.
\end{align}
Then, at points where $H$ resp. $L$ is nonzero,
\begin{align}
\label{rs111}
-\Delta \log H &=2K_1 -2K_2(H-L)\\
\label{rs112}
-\Delta \log L &=2K_1 +2K_2(H-L).
\end{align}
Actually, the Riemann--Hurwitz formula \rf{RH1} can also be deduced from such identities, noting that for instance when $h$ is holomorphic, $L\equiv 0$.
Wolf~\cite{Wolf89} has developed a systematic calculus on the basis of these identities.

The relationship between harmonic and $\pm$ holomorphic maps (or conformal maps, when the target is of higher dimension) is clarified by the following
\begin{lemma}%
\label{l080201}
	Let $\Sigma$ be a Riemann surface, $N$ a Riemannian manifold.
	If $\phi\colon \Sigma \to N$ is harmonic,\index{harmonic} then
	\begin{equation}\label{eq080212}
		T(z) \dz^2 = g\left(\frac{\partial \phi}{\partial z}, \frac{\partial \phi}{\partial z}\right) \dz^2
	\end{equation}
	is a holomorphic quadratic differential.\index{holomorphic quadratic differential}
	Furthermore, $T(z)\dz^2 \equiv 0$ iff $\phi$ is conformal\index{conformal}.
\end{lemma}
The proof of Lemma~\ref{l080201} is an easy computation.
In the case where the target is also a Riemann surface $S$, that is, where we have a map $\phi\colon\Sigma \to S$, we have
\bel{hqd10}
T(z)\dz^2 = \rho^2 \phi_z \overline{\phi}_z\dz^2.
\qe
Then the computation for the proof of Lemma~\ref{l080201} becomes
\begin{align}
	\frac{\partial}{\partial \overline{z}}(\rho^2 \phi_z \overline{\phi}_z)
	=& \rho^2 \overline{\phi}_{z}(\phi_{z\overline{z}} + \frac{2\rho_\phi}{\rho}\phi_z \phi_{\overline{z}})+ \rho^2 \phi_z(\overline{\phi}_{z\overline{z}}+ \frac{2\rho_{\overline{\phi}}}{\rho}\overline{\phi}_z \overline{\phi}_{\overline{z}})\\
	=& 0  \qquad\qquad\text{ by \rf{eq080210b}}.
\end{align}
We obtain, with the above notation $\lambda^2 \dz \dd\overline{z}$ and $\rho^2 \dd\phi \dd\overline{\phi}$
for the metrics on $\Sigma$ and $S$,
\bel{rs113}
\rho^2 \dd\phi \dd\overline{\phi}=T(z)\dz^2+\lambda^2(H+L)\dz\dd\overline{z}+\overline{T}(z)\dd\overline{z}^2.
\qe
(Of course, when we recall the definition \rf{rs110} of $H$ and $L$, the domain metric $\lambda^2 \dz\dd\overline{z}$ drops out.
Also, in our local coordinates,~$T(z)\dz^2=\rho^2 \phi_z \overline{\phi}_{z} \dz^2$.)
Thus, the quadratic differential~$T(z)\dz^2$ is the $(2,0)$-part of the  pullback of the image metric.

In intrinsic terms, a (holomorphic) quadratic differential is a  (holomorphic) section of $T^\ast_{\C} \Sigma \otimes T^\ast_{\C}\Sigma$. Here, $T^\ast_{\C} \Sigma$ is the canonical bundle $K(\Sigma)$ of $\Sigma$.

Every holomorphic quadratic differential\index{holomorphic quadratic differential} on the two-sphere $S^2$
vanishes identically, or putting it differently, the canonical bundle of $S^2$ admits no holomorphic section.
Therefore, Lemma~\ref{l080201} implies
\begin{corollary}\label{c080205}
	Every harmonic map
	\begin{equation}
		\phi \colon S^2 \to N
	\end{equation}
	is conformal, whatever the Riemannian manifold $N$ is.
\end{corollary}
In the next section, we shall see that this corollary reflects the fact that the conformal structure of $S^2$ admits no deformations.

\section{Harmonic maps and Teichm\"uller theory}\label{teich}
We consider Riemann surfaces of genus $p\ge 2$ and the deformations of their complex structure.
For Riemann surfaces of genus 0, there is only one complex structure, that of $\C\mathbb{P}^1$.
For Riemann surfaces of genus 1, tori, there is a complex one-dimensional family of complex structures.
The theory is much simpler than that for genus $p\ge 2$.

As already mentioned, such a Riemann surface $\Sigma$ of genus $p\ge 2$ is a quotient $H/\Gamma$ of the Poincar\'e upper half plane by a discrete group  $\Gamma \in \PSL(2,\mathbb{R})$ of isometries with respect to the hyperbolic metric\index{hyperbolic metric} $\frac{1}{y^2}\dz \otimes \dd\bar{z}$, with $z=x+\ima y$.
Since $\Gamma$ is isomorphic to the fundamental group $\pi_1(\Sigma)$, the Riemann surface is described by a faithful representation $\rho$ of $\pi_1(\Sigma)$ in $\PSL(2,\R)$.
This leads to the approach of Ahlfors and Bers to Teichm\"uller theory.
Representations that only differ by a conjugation with an element of $\PSL(2,\R)$ yield the same conformal structure.
Thus, we consider the space of faithful representations up to conjugacy.
A representation can be defined by the images of the generators, that is, by $2p$ elements of $\PSL(2,\R)$,
and this induces a natural topology on the moduli space\index{moduli space}.
Also, from an easy count, we see that the (real) dimension of the moduli space\index{moduli space} of representations of~$\pi_1(\Sigma)$ in~$\PSL(2,\mathbb{R})$ modulo conjugations is~$6p-6$.
This is the dimension of the moduli space~$\ModuliSpace$ of Riemann surfaces of genus~$p$.
Here, however, we have somewhat more structure than simply the conformal structure of the Riemann surface in question.
In fact, the moduli space presently discussed of discrete, faithful representations of $\pi_1(\Sigma)$ in $\PSL(2,\mathbb{R})$ modulo conjugations yields the Teichm\"uller space~$\Tp$, a simply connected singularity-free infinite cover of $\ModuliSpace$.
Putting it the other way around, $\ModuliSpace$ is obtained as the quotient of $\Tp$ by the mapping class group $\Gamma_p$.
$\Gamma_p$ is the group of homotopy classes of positively oriented diffeomorphisms of the underlying surface $\Sigma$, that is,
\begin{equation}%
\label{mcg1}
	\Gamma_p=\Diff^+(\Sigma)/\Diff_0(\Sigma),
\end{equation}
where $\Diff^+$ stands for orientation preserving diffeomorphisms and $\Diff_0$ for those that are homotopic to the identity of $\Sigma$.

As a moduli of representations of a discrete group,~$\Tp$ also acquires natural structures, like a differentiable and a complex one. It is diffeomorphic, but not biholomorphic to $\C^{3p-3}$.
The complex structure of Teichm\"uller space was originally proposed by Teichm\"uller and investigated by Ahlfors \cite{Ahlfors61a} and Bers \cite{Bers66}, see \cite{Nag88} for details.

The complex dimension $3p-3$ of $\Tp$ equals the dimension of the space of holomorphic quadratic differentials on a surface $\Sigma$ of genus $p$.
In fact, the holomorphic quadratic differentials on $\Sigma$ are the cotangent vectors to $\Tp$ at the point represented by $\Sigma$.

We shall now explore the geometry of $\Tp$ more systematically, with the approach of Riemannian geometry.
Since a two-dimensional Riemannian manifold defines a conformal structure, that is, a Riemann surface, we can naturally look at all Riemannian metrics on a given compact surface $F$ of genus $p$ and then identify those that induce the same conformal structure.
We shall describe here the approach of Tromba and Fischer, see \cite{Tromba92}.
Thus, we consider the space $R_p$ of all smooth Riemannian metrics on $F$.
$R_p$, a space of Riemannian metrics, carries itself a Riemannian metric.
Let $g=g_{ij}\dx^i\otimes\dx^j$ (in local real coordinates~$(x^1,x^2)$) be some Riemannian metric on $F$, that is,~$g$ is an element of~$R_p$.
$(g_{ij})$ is a positive definite symmetric~$2\times 2$ tensor.
Tangent vectors to~$R_p$ at~$g$ then are given by symmetric~$2\times 2$ tensors~$(h_{ij}), (\ell_{ij})$.
The metric of~$R_p$ on the tangent space at~$g$ is then given by
\begin{equation}
\label{eq2.2.6}
	\left( (h_{ij}) , (\ell_{ij})\right)_g \coloneqq\int\limits_F g^{ij} g^{km}h_{ik}\ell_{jm}\sqrt{\det g}\dx^1 \dx^2.
\end{equation}
A word of warning: This $L^2$-metric on the infinite dimensional space $R_p$ is only a weak Riemannian metric, that is, the tangent spaces of $R_p$ are not complete w.r.t.\ this metric. Therefore, the general theory of
Riemann-Hilbert manifolds does not apply. Clarke~\cite{Clarke09b} showed that~$R_p$ becomes a metric space with the distance function induced by this weak metric \rf{eq2.2.6}, and so, the situation is not as bad as one might fear.

In view of our discussion of the gravitino below in Section~\ref{susy2}, we consider \emph{a Riemannian metric as the local difference between the Euclidean coordinate structure and the Riemann surface structure.}
As the uniformization theorem tells us, locally the difference can be made to vanish by choosing suitable coordinates.
Globally, however, the metric on a compact Riemann surface cannot be made Euclidean, unless the surface is a torus, i.e., of genus one.
Nevertheless, the metric carries some redundant information about the underlying conformal structure, as we shall now analyze.

Recall that the Riemannian manifolds \((F, g_1)\) and \((F, g_2)\) are called isometric if there is a diffeomorphism \(k\colon F\to F\), such that \(k^*g_2=g_1\).
In that case we would like to consider the metrics \(g_1\) and \(g_2\) as equivalent.
Conversely, any diffeomorphism \(k\colon F\to F\) yields the obviously isometric Riemannian manifolds~\((F, g)\) and \((F, k^*g)\).
This is not yet reflected in the definition of $R_p$.
We therefore divide out of \(R_p\) the action of the (orientation preserving) diffeomorphism group $D_p$ of $F$.
The diffeomorphism group \(D_p\) acts isometrically on \(R_p\) and hence the Riemannian metric~\eqref{eq2.2.6} descends to the quotient of \(R_p\) by \(D_p\).

However, even metrics that are not isometric might still induce the same conformal structure, and we therefore should also identify all such metrics.
Such metrics differ by a some scalar factor.
When we multiply a metric $g$ by a positive function~$\lambda$, the metric~$\lambda g$ induces the same conformal structure as $g$, and  conversely.

In order to eliminate  the ambiguity of the conformal factor, one seeks a slice in $R_p$ transversal to the conformal changes.
As explained above, by the  Poincaré uniformization theorem, any Riemannian metric on a surface of genus $p>1$ is conformally equivalent to a unique hyperbolic metric, that is, one with curvature $-1$, turning it into a quotient $H/\Gamma$.
The moduli space $\ModuliSpace$ is then obtained as the space $R_{p,-1}$ of metrics of curvature \(-1\) divided by the action of the diffeomorphism group $D_p$.
The geometric structures on $R_p$ then induce corresponding geometric
structures on $\ModuliSpace$, see~\cite{Tromba92}.
We can use the metric~\rf{eq2.2.6} to identify those directions on $R_p$ that are orthogonal to the conformal slice and the action of the diffeomorphism group.
These are the ones that lead to nontrivial deformations of the complex structure.
We recall that a tangent vector to $R_p$ is a symmetric $2\times 2$ tensor $(h_{ij})$.
Such a tensor is orthogonal to conformal multiplications when trace-free, and orthogonal to the action of $D_p$ when divergence-free.
Trace- and divergence-free symmetric two-tensors are nothing but holomorphic quadratic differentials on the Riemann surface, and these are precisely the objects produced from harmonic mappings, see Lemma~\ref{l080201}.

The Riemannian metric on $R_p$ then induces a Riemannian metric on the Teichm\"uller space~$\Tp$.
This theory is described in detail in~\cite{Tromba92}.
And it follows from the preceding considerations that this metric in turn induces a product between holomorphic quadratic differentials on the Riemann surface~$\Sigma$ in question.
In complex notation, with the hyperbolic metric on~$\Sigma$ denoted by~$\lambda^2\dz\dd\overline{z}$, let $q_1,q_2$ be holomorphic quadratic differentials on~$\Sigma$.
Their product w.r.t.\ the metric then is given by
\begin{equation}%
\label{1.2.1}
	{(q_1,q_2)}_g = 2 \operatorname{Re}\int q_1\bar{q}_2  \frac{1}{\lambda^2(z)}\frac{i}{2}\dz \wedge \dd\overline{z}.
\end{equation}
This is the Weil--Petersson metric $g_{WP}$.
The Weil--Petersson is a K\"ahler metric w.r.t.\ the complex structure of the moduli space.
This metric can be studied with the help of harmonic mappings.
For instance, Tromba~\cite{Tromba96,Tromba96a} found that the energy functional for harmonic maps between Riemann surfaces as a function of the conformal structure of the domain yields a strictly convex exhaustion function for the Weil--Petersson metric on $\Tp$.
For some reviews of the geometry of the Weil--Petersson metric, see~\cite{Wolpert05,Wolpert08}.

Tromba~\cite{Tromba86} proved that $g_{WP}$ has negative sectional curvature, and its holomorphic sectional curvature even has a negative upper bound $k<0$, and this can again be derived with the help of harmonic mappings, see~\cite{Jost91a,Jost/Peng92}.

We now explain the harmonic map approach to Teichm\"uller theory in some more detail.
We consider harmonic maps~$\phi\colon\Sigma \to S$ between compact Riemann surfaces of the same genus $\ge 2$. In contrast to the domain~$\Sigma$, for which we only need a conformal structure, the target~$S$ needs a Riemannian metric, and we can take the hyperbolic metric of constant curvature $-1$.
As explained, we equip the surfaces with an additional structure, a representation of the fundamental group. Equivalently, we fix homotopy classes of diffeomorphisms. We then assume that our harmonic map respects these homotopy classes. In still other words, both~$\Sigma$ and~$S$ are considered to be diffeomorphic to some abstract underlying differentiable surface~$F$, and we assume that the harmonic map~$\phi$ is then homotopic to the identity of~$F$.

By Lemma~\ref{l080201}, the harmonic map $\phi$ induces a holomorphic quadratic differential on $\Sigma$, that is, a holomorphic section of~$T^\ast_{\C} \Sigma \otimes T^\ast_{\C}\Sigma$. So, the harmonic map~$\phi$ yields a cotangent vector of the Teichm\"uller space~$\Tp$ at the point corresponding to~$\Sigma$.
And this cotangent vector then varies when we vary the conformal structure of either $\Sigma$ or $S$. Let us consider the effect of variations of the latter. From Theorem~\ref{diffeo} and Lemma~\ref{l080201} and basic harmonic map theory, we have
\begin{lemma}
Let $\Sigma$ be a compact Riemann surface with local conformal coordinate $z=x+\ima y$, $S$ another such surface of the same genus $p\ge 2$, equipped with its unique hyperbolic metric, locally written as $\rho^2(h)\dd h \dd\overline{h}$. We also fix a homotopy class of (orientation preserving) diffeomorphisms $k\colon\Sigma \to S$. Then there exists a unique harmonic diffeomorphism
\begin{equation}
\label{harm50}
	\phi =\phi (\Sigma,S)\colon\Sigma \to S
\end{equation}
homotopic to $k$.
\begin{equation}\label{harm51}
	T(z)\dz^2 = \rho^2 \phi_z \overline{\phi}_{z} \dz^2
\end{equation}
is a holomorphic quadratic differential, and $\phi$ is conformal iff $T=0$.
\end{lemma}
We keep $\Sigma$ fixed. The harmonic map $\phi$ then varies with~$S$, and so then does the holomorphic quadratic differential on~$\Sigma$. We obtain a map
\bel{harm52}
q(\Sigma)\colon \Tp \to Q(\Sigma)
\qe
into the space of holomorphic quadratic differentials~$Q(\Sigma)$ on~$\Sigma$.
The map~$q$ is defined on Teichm\"uller space~$\Tp$ instead of on the
moduli space~$\ModuliSpace$ because the harmonic map depends on the choice of a
diffeomorphism $k\colon \Sigma \to S$. Wolf \cite{Wolf89} (see also \cite{Jost91}) then showed
\begin{theorem}%
\label{Wolf}
For any $\Sigma$, $q(\Sigma)$ is a bijection between $\Tp$ and $Q(\Sigma)$.
\end{theorem}
The proof makes crucial use of the identities \rf{rs111}, \rf{rs112}.
Theorem~\ref{Wolf} says that with its natural differentiable structure, $\Tp$ is diffeomorphic to the vector space $\C^{3p-3}$.
This is Teichm\"uller's theorem.
As explained above, however, it is not biholomorphic to that space with respect to its natural complex structure.
Also, $q(\Sigma)$ is not an isometry w.r.t.\ the Weil--Petersson metric on $\Tp$.
Nevertheless, there does exist a relationship between harmonic maps and the Weil--Petersson metric that we are now going to describe.
For that purpose, we consider the energy
\bel{harm53}
E(\phi(\Sigma,S))=\int_\Sigma \rho^2 (\phi(z))(\phi_z
\overline{\phi}_{\overline{z}}+\phi_{\overline{z}}\overline{\phi}_z)\frac{i}{2}\dz \wedge \dd\overline{z}=\int_\Sigma (H+L)\lambda^2(z)\frac{i}{2}\dz \wedge \dd\overline{z}
\qe
as a function of $S$ and present some formulas obtained in \cite{Wolf89,Jost91}.
When $S=\Sigma$ as an element of $\Tp$, that is, when the conformal structure of $S$ is such that $k$ is homotopic to a conformal diffeomorphism, and then the harmonic map $\phi$ is conformal.
Therefore, $L=0$, and the integral of the Jacobian,
\bel{harm54}
J(h)=\int_\Sigma (H-L)\lambda^2(z)\frac{i}{2}\dz \wedge \dd\overline{z}
\qe
is a topological invariant that does not depend on the conformal structures of $\Sigma$ and $S$, but only on their genus. Therefore, comparing with \rf{harm53}, $E$ as a function of $S$ achieves its
minimum precisely when $\phi$ is conformal, that is, when the conformal structures of $\Sigma$ and $S$ are the same. At the point $S=\Sigma \in \Tp$, the infinitesimal variations of the target surface
correspond to the cotangent vectors of $\Tp$ at $\Sigma$, that is, to the holomorphic quadratic differentials on $\Sigma$. Since $E$ achieves its minimum here,
\bel{harm55}
E_q =0 \text{ for all } q \in Q(\Sigma),
\qe
where $E_q = \left.\frac{d}{dt}\right|_{t=0} E(\phi(\Sigma, S+tq))$ is the derivative of $E$ for a variation of the target in the direction of $q$, and where we recall the identification of $\Tp$ with~$Q(\Sigma)$ given by Theorem~\ref{Wolf}.
Moreover,
\begin{lemma}
	$S=\Sigma$ is the only critical point of $E$ as a function of $S\in \Tp$.
	$E(\phi(\Sigma,\cdot))$ is a proper function on $\Tp$.
	Its second derivatives at $S=\Sigma$ for $q_1,q_2 \in Q(\Sigma)$ are
	\bel{harm56}
	E_{q_1 \overline{q}_2}= 2\mbox{Re} \;\int_\Sigma q_1 \overline{q}_2 \frac{1}{\lambda^2}\frac{i}{2}\dz \wedge \dd\overline{z}.
	\qe
	Consequently, $E(\phi(\Sigma,\cdot))$ is a proper exhaustion function with a single critical point at which its Hessian is positive definite.
\end{lemma}
Equation~\rf{harm56} tells us  that the second variation of the energy at $S=\Sigma$ yields the Weil--Petersson product of the holomorphic quadratic differentials $q_1,q_2$.
On this basis, one can develop systematic expansions~\cite{Jost91a}.

\section{Variations of conformal structures}\label{vari}

We now look at the case where we vary the domain instead of the target. The corresponding calculations extend to the case where we have a general Riemannian manifold $N$ as a target. This is useful when one studies minimal surfaces in Riemannian manifolds, see \cite{Jost91a,Jost97c}. For simplicity, we assume here that all the harmonic maps into $N$ that we consider here are unique in their homotopy classes, so that we don't have to deal with issues of bifurcation.

There is a new difficulty, however. When we vary the conformal structure of the domain, then $z=x+\ima y$ no longer remain holomorphic coordinates. Thus, when we wish to express our formulas in the elegant complex notation, we need to account for that. Thus, instead of expressing things in non-holomorphic coordinates, we pull everything back from the varying surfaces to the fixed reference surface that we want to vary. And for our computations, we shall work with a variation of the domain metric instead of the conformal structure. While, as explained, many of the variations of the metric, more precisely those by diffeomorphisms or conformal factors,  are redundant, insofar as they don't change the underlying complex structure, this is computationally easier. Thus, let $\gamma=\lambda^2 \dz\otimes\dd\overline{z}$ be a conformal metric on our domain $\Sigma$. And let $\gamma_t, t\in (-\varepsilon, \varepsilon)$ be a smooth variation of it. By what we have said before, we might assume that all the metrics $\gamma_t$ are hyperbolic metrics. Thus, the surfaces $\Sigma_t$ with the conformal structure determined by $\gamma_t$ (via Gauss' uniformization theorem) can be represented as quotients $H/\Gamma_t$ of the upper half plane. While this can help with the computations, here we only display the results, referring to \cite{Jost91a}, Sections 3.3 and 6.4, for details.

For a mapping $u\colon\Sigma \to N$ and a metric $\eta$ on $\Sigma$, we let $E(u,\eta)$ be its energy~\rf{eq080211}. Again, to repeat an important point, although in the second line of \rf{eq080211}, the metric no longer appears explicitly, it is still implicitly contained, because the holomorphic coordinates $z=x+\ima y$ depend on the metric.

Also, let $u(\gamma_t)\colon \Sigma_t \to N$ be the harmonic maps that vary the harmonic map $u(\gamma)\colon\Sigma \to N$ (it follows from general principles that these maps depend smoothly on $t$).
We let $\kappa_t\colon \Sigma \to \Sigma_t$ be diffeomorphisms, again depending smoothly on $t$, with $\kappa_0 =\mathrm{id}$.
Since $u=u(\gamma)$ is harmonic, we consider $u\circ \kappa_t^{-1}\colon \Sigma_t \to N$ as an approximation of the harmonic $u(\gamma_t)$.
Since, as $u_0$ is harmonic,
\bel{vari1}
\frac{\partial E(u_t,\gamma_t)}{\partial u}=0 \text{ at } t=0
\qe
for any smooth variation $u_t$ of $u=u_0$, this won't make a difference  for the first order computations to follow.  Indeed, with \rf{vari1}, we compute
\bel{vari2}
\frac{d}{dt} E(u\circ \kappa_t^{-1}, \gamma_t)_{|t=0}=-2\mbox{Re} \int T(z) \mu_{\overline{z}}\frac{i}{2}\dz \wedge \dd\overline{z}
\qe
with
\begin{equation}%
\label{vari3}
	T(z)\dz^2 = g\left(\frac{\partial \phi}{\partial z},\frac{\partial\phi}{\partial z}\right) \dz^2
\end{equation}
being the holomorphic quadratic differential associated with the harmonic map~$u$, see~\rf{eq080212}, and with
\bel{vari4}
\mu =\frac{d \kappa_t}{dt} \text{ at } t=0
\qe
being the infinitesimal variation of the underlying structure. Now,
\bel{vari5}
Q(z)\dz^2\coloneqq \lambda^2 \overline{\mu_{\overline{z}}}\dz^2
\qe
transforms as a quadratic differential, and we may write \rf{vari2} as
\bel{vari6}
\frac{d}{dt} E(u\circ \kappa_t^{-1}, \gamma_t)_{|t=0}=-2\mbox{Re} \int T(z)\overline{Q}(z)\frac{1}{\lambda^2} \frac{i}{2}\dz \wedge \dd\overline{z}.
\qe
This is the negative Weil--Petersson product~\rf{1.2.1} between the quadratic differentials~$T(z)\dz^2$ and $Q(z)\dz^2$. Here, $T\dz^2$ is holomorphic when $u$ is harmonic, but $Q$ needs not be holomorphic. As we have discussed above, however, those variations $\mu$ in \rf{vari4} that do not change the conformal structure lead to quadratic differentials $Q\dz^2$ in \rf{vari5} that are orthogonal to the holomorphic ones, and therefore, their Weil--Petersson products \rf{vari6} with a holomorphic $T\dz^2$ vanish.

Conceptually, we can also consider the left hand side of \rf{vari6} as the variation of $\gamma$ in the direction $Q\dz^2$ and write it as
\bel{vari7}
\frac{d}{d\gamma} E(u(\gamma),\gamma)(Q\dz^2),
\qe
recalling again that for a harmonic $u$, we have $\frac{\partial}{\partial u} E(u,\gamma)=0$. We can then also proceed to compute second derivatives when $u(\gamma)$ is again assumed to be harmonic, and in fact, we get the analogue of \rf{harm56},
\bel{vari8}
\frac{d^2}{d\gamma^2} E(u(\gamma),\gamma)(Q_1\dz^2,Q_2\dz^2)= 2\mbox{Re} \int Q_1(z)\overline{Q_2}(z)\frac{1}{\lambda^2} \frac{i}{2}\dz \wedge \dd\overline{z},
\qe
see \cite{Jost91a}, Section 6.4.

\section{Dirac-harmonic maps and the gravitino}
From the perspective of quantum field theory, the harmonic map and Riemann surface model is only the shadow of a much deeper and richer one, the supersymmetric nonlinear $\sigma$-model.
The latter, while not directly a model of real physical forces, is one of the most basic and important models of quantum field theory, because of the symmetries it encodes and the conceptual and computational links with other models it provides.
And even more fundamentally, it has become the basic action functional of superstring theory.
For mathematical introductions, see~\cite{Deligne99a,Deligne99b,Jost09a}.

As we have seen, the fundamental symmetry, conformal invariance, only holds when the domain is two-dimensional, that is, a Riemann surface.
Another symmetry, supersymmetry, requires an additional field besides the harmonic map, a twisted spinor field along the map.
It also requires that this field be anticommuting.
This requires the introduction of additional anticommuting, so-called odd variables, in addition to the standard coordinate functions that we have been working with so far.
This brings us into the realm of supergeometry, which has a number of features that appear odd from the usual perspective of differential geometry.
In particular, for anticommuting variables, one cannot form inequalities, and therefore, the standard tools and estimates of PDE theory no longer apply.
In order to prepare the ground, however, we shall first  develop a mathematical version of the theory that need only ordinary, commuting variables and fields, introduced in~\cite{Chen/Jost/Li/Wang05b,Chen/Jost/Li/Wang05a}.
This has since become a very active direction of research in the calculus of variations and geometric analysis, with many difficult and challenging problems that lead to the development of new and powerful mathematical tools that are quite useful also for other problems in geometric analysis.
We now want to describe that theory.

\subsection{Dirac-harmonic maps}
$(N,g)$ still is a Riemannian manifold.
Let $(\Sigma,\gamma)$ be a Riemann surface with a metric $\gamma$, as before, but now also equipped with a spin structure. Metrics and connections on the bundles appearing in the sequel will be those induced by~$\gamma$ and its Levi-Civita connection.
The~$\SpinGroup$ principal bundle corresponding to the spin structure is called $\SpinStructure$.
$\CliffordBundleM$  is the corresponding Clifford algebra bundle; it is  isomorphic to the quotient of the tensor algebra by the two-sided ideal generated by
\begin{equation}
	X\otimes Y + Y\otimes X + 2\gamma(X, Y),
\end{equation}
for $X, Y\in\Gamma(T\Sigma)$.
Its  fiber is the Clifford algebra $\CliffordAlgebraM$.
We consider the spinor bundle $\SpinorBundleM=\SpinStructure\times_{\SpinGroup}\CliffordAlgebraM$ where $\SpinGroup\subset\CliffordAlgebraM$ acts via the left-regular representation of $\CliffordAlgebraM$.
We denote the Clifford multiplication of a vector \(X\) with a spinor \(s\) by \(X\cdot s\).
For more details on the spin geometry and the spinor bundle, we refer to~\cite{Jostetal17}.
General background material can be found in~\cite{Jost17}.

We shall study the action functional $\ADH$ defined on the space
\begin{equation}
	\mathcal{X}(\Sigma,N)=\{(\phi,\psi)\big| \phi\colon\Sigma\to N, \psi\in\Gamma(\SpinorBundleM\otimes \phi^*TN)\},
\end{equation}
where $\Gamma(\SpinorBundleM\otimes \phi^*TN)$ is the space of sections of the twisted spinor bundle \(\SpinorBundleM\otimes \phi^*TN\).
The action functional then is
\begin{equation}
\label{A-intro}
	\ADH(\phi, \psi;\gamma, \chi) \coloneqq \int_\Sigma ( \|d\phi\|^2 + \langle \psi, \Dirac \psi \rangle_{\SpinorBundleM\otimes \phi^*TN}) \dv.
\end{equation}
Here, the first term is the same as in \rf{eq080106}.
In the second term, $\Dirac$ is the  twisted spin Dirac operator on $\SpinorBundleM\otimes \phi^*TN$.
In local coordinates $y^j$ on $N$, \(\psi=\psi^j\otimes\phi^*\frac{\partial}{\partial y^j}\) and with a local orthonormal frame $e_\alpha$ on $\Sigma$, it is given by
\begin{equation}%
\label{dirac}
	\Dirac\psi = e_\alpha \cdot \nabla^{\SpinorBundleP\otimes\phi^*TN}_{e_\alpha}\psi
	= \dirac \psi^j \otimes \phi^*(\frac{\partial }{\partial y^j}) + e_\alpha \cdot \psi^j \otimes \nabla^{\phi^* TN}_{e_\alpha} \phi^*\frac{\partial}{\partial y^j}
\end{equation}
where $\dirac$ is the ordinary spin Dirac operator on \(\SpinorBundleM\) (see for instance~\cite{Jost17}).
Furthermore, here and in the following we use the convention that we sum over repeated lower indices of the orthogonal frames \(e_\alpha\).
Since we choose the Clifford algebra $\CliffordAlgebraM$
\begin{equation}%
\label{intpar}
	\int \langle \psi, \Dirac \psi \rangle_{\SpinorBundleM\otimes \phi^*TN} \dv =\int \langle \Dirac \psi, \psi \rangle_{\SpinorBundleM\otimes \phi^*TN} \dv
\end{equation}
under integration by parts as well as by symmetry of the product $\langle \cdot,\cdot \rangle_{\SpinorBundleM\otimes  \phi^*TN}$, that is, the corresponding term in the action functional is symmetric.
Had we chosen $\CliffordAlgebraP$, integration by parts would have introduced a minus sign,  and the corresponding term in the action functional would have been zero.

From \rf{dirac}, we see that $\Dirac$  depends not only on the spin geometry of the domain $\Sigma$, but also on the map $\phi$.
This leads to a coupling of the two fields $\phi$ and $\psi$ in \rf{A-intro}.
Therefore, also the Euler--Lagrange equations for \rf{A-intro} are coupled.
They are
\begin{align}%
\label{EL1}
	\tau(\phi) &= \frac{1}{2}\mathcal{R}^N\left(\phi, \psi\right), \\%
\label{EL2}
	\Dirac \psi &= 0.
\end{align}
Here, \(\tau(\phi)\) is the tension field~\eqref{ha30} of \(\phi\), the Dirac operator~\(\Dirac\) was defined in~\eqref{dirac}, and \(\mathcal{R}^N\left(\phi, \psi\right)\) is a term depending on the Riemannian curvature%
\footnote{%
	We use the sign conventions of~\cite{Jost17}, i.e.,
	\begin{equation}
		\Rm_{ijkl} = g\left(\Rm\left(\frac{\partial}{\partial y^k}, \frac{\partial}{\partial y^l}\right) \frac{\partial}{\partial y^j}, \frac{\partial}{\partial y^i}\right)
		= g\left(\nabla_{\frac{\partial}{\partial y^k}}\nabla_{\frac{\partial}{\partial y^l}}\frac{\partial}{\partial y^j} - \nabla_{\frac{\partial}{\partial y^l}}\nabla_{\frac{\partial}{\partial y^k}}\frac{\partial}{\partial y^j}, \frac{\partial}{\partial y^i}\right).
	\end{equation}
} on \(N\):
\begin{equation}
	\mathcal{R}^N\left(\phi, \psi\right) = \left<\psi^i, e_\alpha\cdot\psi^j\right>_\SpinorBundleM R^{\phi^*TN}\left(\phi^*\frac{\partial}{\partial y^i}, \phi^*\frac{\partial}{\partial y^j}\right) \phi_* e_\alpha.
\end{equation}
These equations were derived in~\cite{Chen/Jost/Li/Wang05b}, and the investigation of the regularity of their solutions has been started in~\cite{Chen/Jost/Li/Wang05a}.
Further aspects and extensions of the model have been studied in~\cite{Branding15,Branding16b,Chen/Jost/Sun14,Chen/Jost/Sun/Zhu15,Chen/Jost/Sun/Zhu17b,Chen/Jost/Wang07a, Chen/Jost/Wang08b,Chen/Jost/Wang/Zhu11,Jost/Liu/Zhu16,Jost/Liu/Zhu17c,Jost/Liu/Zhu17e,Sharp/Zhu16,Wang/Xu09,Zhu09,Z-DHMDSSNSC} and many other papers.
For the regularity theory, the method of Rivi\`ere~\cite{Riviere07} turned out to be very useful.

The existence theory met with the difficulty that, in contrast to the energy~\rf{eq080106}, the action \rf{A-intro} is not bounded from below.
This comes from the term~$\langle \psi, \Dirac \psi \rangle$, and ultimately from the fact that the spectrum of the Dirac operator~$\dirac$, in contrast to that of the Laplace operator \rf{eq020112}, is not bounded from below.
Therefore, variational methods do not apply to show the existence of solutions.
As already mentioned at the end of Section \ref{harmonic}, an alternative to variational methods consists in heat flow methods.
In the present case, however, there is the additional difficulty that the Dirac operator as a first order differential operator does not admit a natural parabolic version, again in contrast to the Laplace operator.
Therefore, in \cite{Chen/Jost/Sun/Zhu17}, a novel elliptic-parabolic system has been introduced.
The idea is to convert the second order elliptic system \rf{EL1} into a parabolic system, but to carry the first order elliptic system \rf{EL2} along as an elliptic side constraint.
This leads to the system
\begin{align}%
\label{EL3}
	\frac{\partial \phi}{\partial t}=&\tau(\phi) - \frac{1}{2}\mathcal{R}^N\left(\phi, \psi\right), \\%
\label{EL4}
	\Dirac \psi =& 0,
\end{align}
with suitable initial conditions, and also suitable boundary conditions when the domain has a boundary, see
\cite{Chen/Jost/Wang/Zhu11}. The short and long time existence of solutions, their regularity and their asymptotic behavior has been further investigated in \cite{Jost/Liu/Zhu17a,Jost/Liu/Zhu17b,Jost/Liu/Zhu17d,Wittmann17}.

\subsection{The gravitino}\label{gravi1}
In fact, the full supersymmetric action functional (which we shall discuss in Section~\ref{susy2} below) has a still richer structure and involves additional fields, in particular the gravitino field, a superpartner of the domain metric $\gamma$.
Recently, in  \cite{Jostetal17},  we have constructed and investigated a version of this model that involves only commuting variables, analogous to and generalizing the Dirac-harmonic action function \rf{A-intro}.
For some further results on this model, see~\cite{Jost/Wu/Zhu17a,Jost/Wu/Zhu17b,Jostetal17b}.

We shall now describe that model.
\begin{defi}
	A \emph{gravitino} $\chi$ is a section of the bundle \(\SpinorBundleM\otimes TM\).
\end{defi}
For the moment, we assume $\chi$ to be smooth, but smoothness will become an issue below. In \cite{Jost/Kessler/Tolksdorf16}, a gravitino is introduced as a section of the bundle
$\SpinorBundleM\otimes T^*M$; we can, however, use the metric~$\gamma$ to identify~$T^*M$ with~$TM$.

The bundle $\SpinorBundleM\otimes TM$ can be decomposed into irreducible representations of $\SpinGroup$. In fact, we get  two representations of type $\frac{1}{2}$ and two of type $\frac{3}{2}$.
We define the operator $Q$ as the projection onto the $\frac{3}{2}$-parts. In local coordinates, it is given by
\bel{Qdef}
	Q\chi = -\frac{1}{2} e_\alpha \cdot e_\beta \cdot \chi^\alpha \otimes e_\beta.
\qe

We can now present the action functional, whose form, however, will only be derived and justified and whose content will only become clear in the next Section \ref{susy2},
\begin{equation}
\label{eq:AF}
	\begin{split}
		\ADHG(\phi, \psi;\gamma,\chi)&\coloneqq \int_\Sigma \left(\vphantom{\frac16}|\dd \phi|_{T^*\Sigma\otimes \phi^*TN}^2 + \langle \psi, \Dirac \psi \rangle_{\SpinorBundleM\otimes \phi^*TN} \right.\\
			&\hspace{4em} -4\langle (\mathrm{Id}\otimes\phi_*)(Q\chi), \psi \rangle_{\SpinorBundleM\otimes \phi^*TN} \\
			&\hspace{4em}\left.- |Q\chi|^2_{S\otimes TM} |\psi|^2_{\SpinorBundleM\otimes \phi^*TN}
			-\frac{1}{6} \Rm^{N}(\psi)\right) \dv_\gamma,
	\end{split}
\end{equation}
with the curvature term%
\begin{equation}
	\Rm^{N}(\psi)
	=\Rm^{\phi^*TN}_{ijkl}\langle \psi^i, \psi^k \rangle_S \langle \psi^j,\psi^l \rangle_S
	=\langle SR(\psi),\psi\rangle_{\SpinorBundleM\otimes \phi^*TN}
\end{equation}
where
\begin{equation}%
\label{def-SR}
	SR(\psi)\coloneqq\langle \psi^l, \psi^j \rangle_S \psi^k
	\otimes \phi^*(\Rm^{\phi^*TN}(\frac{\p}{\p y^k},\frac{\p}{\p y^l})\frac{\p}{\p y^j}).
\end{equation}
Since  $Q$ is an orthogonal projection,
\bel{Qproj} |Q\chi|^2_{S\otimes T\Sigma}= \langle \chi, Q\chi \rangle.
\qe
And we observe that in the action functional \rf{eq:AF}, only the $\frac{3}{2}$-part $Q\chi$ of $\chi$ enters.

The functional $\ADHG(\phi,\psi;\gamma,\chi)$ is of geometric interest because of its symmetries.
Like the harmonic map energy \rf{eq080106},  which it generalizes, it is invariant under generalized conformal transformations of the metric.
However, as the spinor bundles depend on the metric we need to be a little careful.
It is explained in~\cite{BG-SODVM} that for any two metrics \(\gamma\) and \(\gamma'\) there is an isometric isometry \(b\colon (T\Sigma, \gamma)\to (T\Sigma, \gamma')\) that lifts to an isometry \(\beta\colon \SpinorBundleM\to\SpinorBundleM'\) of the spinor bundles.
Using those maps, we obtain for \(\gamma'=e^{2u}\gamma\)
\begin{equation}
	\ADHG(\phi,e^{-\frac{1}{2}u}(\beta\otimes\mathrm{Id})\psi; e^{2u}\gamma, e^{-\frac{1}{2}u}(\beta\otimes b)\chi)
	=\ADHG(\phi,\psi;\gamma,\chi).
\end{equation}
Here we use the scaling behavior of the twisted Dirac operator, $\Dirac^{e^{2u}\gamma}e^{-u}\psi = e^{-2u}\Dirac^\gamma\psi$.
However, in our setting of conformal rescaling of the metric, the maps can be realized by \(b=e^{-u}\mathrm{Id}\) and \(\beta=e^{-\frac{1}{2}u}\mathrm{Id}\) which leads to the following more explicit invariance formula:
\begin{equation}
	\ADHG(\phi,e^{-u}\psi;e^{2u}\gamma,e^{-2u}\chi) =\ADHG(\phi,\psi;\gamma,\chi).
\end{equation}
We call this invariance rescaled conformal invariance, since not only the metric but also the spinors are rescaled.

A further invariance comes from the super Weyl transformations, whose role will again become clear in the next Section \ref{susy2},
\begin{equation}
	\ADHG(\phi,\psi;\gamma,\chi+\chi')=\ADHG(\phi,\psi;\gamma,\chi)
\end{equation}
with $Q\chi'=0$. As we have noted above, only the part $Q\chi$ enters into $\ADHG$, and since $Q$ is a projection, we have $Q^2\chi =Q\chi$.

As for the energy functional, we also have diffeomorphism invariance, again to be generalized in the next Section \ref{susy2}, and furthermore, $\ADHG$ is \(\SpinGroup\)-gauge-invariant.

Although, as we shall see in Section~\ref{susy2}, the action functional comes from that of the two-dimensional nonlinear supersymmetric sigma model, we are using here  commuting spinors, and therefore~\eqref{eq:AF} is in general no longer supersymmetric.
Only in the case \(\chi=0\) and under certain conditions on the curvature of the target, the functional~\(\ADHG\) is supersymmetric also in the case of commuting spinors, see~\cite{Jostetal17b}.
Thus, we have to pay a price for moving from anticommuting to commuting variables.
That price is a loss of symmetry.
Nevertheless, as described, our functional still possesses many symmetries, and these are crucial for its analysis and for its geometric content.
Furthermore, as described in~\cite{KT-FSRSSCS}, the functional \(\ADHG\) is completely determined by the requirement of rescaled conformal and super Weyl invariance, given that the equations of motion are at most of second order.

From a somewhat lengthy computation, in \cite{Jostetal17}, and adopting the notation there, we have obtained
\begin{theo}
The system of Euler--Lagrange equations for the action functional~$\ADHG$ is
\begin{equation}%
\label{EL-eq}
	\begin{split}
		\tau(\phi) &= \frac{1}{2}\mathcal{R}^N\left(\phi, \psi\right) - \frac{1}{12}S\na R(\psi)   \\
			&\quad -(\langle \na^S_{e_\beta}(e_\al \cdot e_\beta \cdot \chi^\al), \psi \rangle_S
			+ \langle e_\al \cdot e_\beta \cdot \chi^\al, \na^{\SpinorBundleM\otimes \phi^*TN}_{e_\beta} \psi \rangle_S),  \\
		\Dirac\psi &= |Q\chi|^2\psi +\frac{1}{3}SR(\psi)+2(\mathrm{Id}\otimes \phi_*)Q\chi.
	\end{split}
\end{equation}
\end{theo}
These equations already make the growth order transparent with which the various fields enter.
The term $SR(\psi)$, which  has been defined in \rf{def-SR},  is cubic in~$\psi$;  $S\na R(\psi)$ involves derivatives of the curvature tensor $R^N$ and is quartic in~$\psi$.

As for the harmonic map system, one can define weak solutions. Here, as before, $\phi$ needs to be of class $W^{1,2}$, and the spinor $\psi$ needs to be of class $W^{1,4/3}$ (this can be seen from the action functional and the Sobolev embedding theorem; the latter guarantees that such a $\psi$ is, in particular, of class $L^4$).
In~\cite{Jostetal17}, we have shown
\begin{theo}%
\label{theorem 1}
Assume that the metric $\gamma$ and the gravitino $\chi$ are smooth. Then
	the critical points of the  action functional $\ADHG$
	are smooth.
\end{theo}
Compare the discussion at the end of Section \ref{func}. Here, however, the regularity theory is much harder. A crucial ingredient is Rivi\`ere's method \cite{Riviere07}.

As before, the symmetries lead to conserved Noether currents, as we have already seen that conformal invariance of the harmonic map energy leads to a holomorphicity of the quadratic differential.
Now in our situation here we have more symmetries, hence we may expect more conserved currents. They are actually given by the variation of the action functional with respect to the Riemannian metric and the gravitino respectively. Formally, we have
\begin{align}
	\partial_\gamma \ADHG(\phi,\psi;\gamma,\chi)(\delta \gamma) &=\int_{\Sigma} \left<\delta\gamma, T\right> \dv_\gamma\\
	\partial_\chi \ADHG(\phi,\psi;\gamma,\chi)(\delta \chi) &= \int_{\Sigma} \left<\delta\chi, J\right> \dv_\gamma.
\end{align}
In a local orthonormal frame $(e_\alpha)$, the \emph{energy-momentum tensor} $T$ is given by $T=T_{\alpha\beta}e^\alpha\otimes e^\beta$ with components
\begin{align}
	T_{\alpha\beta}=
	&2\langle\phi_*e_\alpha,\phi_*e_\beta\rangle+\frac{1}{2}\left\langle\psi,e_\alpha\cdot\nabla^{S\otimes\phi^*TN}_{e_\beta}\psi + e_\beta\cdot\nabla^{S\otimes\phi^*TN}_{e_\alpha}\psi\right\rangle \\
	&+\langle e_\eta\cdot e_\alpha\cdot\chi^\eta\otimes\phi_* e_\beta +e_\eta\cdot e_\beta\cdot\chi^\eta\otimes\phi_* e_\alpha,\psi\rangle \\
	&-\left(|\dd\phi|^2+\langle\psi,\Dirac\psi\rangle-4\langle(\mathrm{Id}\otimes\phi_*)Q\chi,\psi\rangle-|Q\chi|^2|\psi|^2-\frac{1}{6}\Rm(\psi)\right)\gamma_{\alpha\beta}
\end{align}
while the \emph{supercurrent} $J$ is given by $J=J^\alpha\otimes e_\alpha$ with components
\begin{equation}
	J^\alpha=2\langle\phi_* e_\beta, e_\beta\cdot e_\alpha\psi\rangle+|\psi|^2e_\beta\cdot e_\alpha\cdot\chi^\beta.
\end{equation}
The rescaled conformal symmetry prescribes the $\gamma$-trace of the energy-momentum tensor $\trace(T)$.
However, since the rescaled conformal invariance is not a pure conformal invariance, the trace will not be zero but rather depend on \(\psi\) and \(\chi\).
The super Weyl symmetry prescribes the vanishing of the $\frac{1}{2}$-part of the supercurrent: $(\mathrm{Id}-Q)J=0$.
As in the case of harmonic maps, the diffeomorphism invariance yields a conservation law of divergence-type:
\begin{equation}
	\diverg_\gamma(T)+\diverg_\chi(J)=0.
\end{equation}
Here $\diverg_\chi$ is defined as the formal adjoint operator of the Lie derivative operator $L^{S\otimes TM}$, just as $\diverg_\gamma$ is formally adjoint to the Lie derivative $L$ on symmetric two tensors on $\Sigma$, see \cite{Tromba92, Jostetal17b}.
This is a law which involves the derivatives of the currents.

As explained above, we cannot expect full supersymmetry in this model.
However, in the rare cases, that the action~\(\ADHG\) is supersymmetric, we obtain another conservation law of divergence type.
In that case the energy-momentum tensor can be considered as a holomorphic quadratic differential and the metric dual of the supercurrent as a holomorphic section of the \(\frac32\)-part of \(\SpinorBundleM^*\otimes T^*\Sigma\).
This issue is analyzed in detail in \cite{Jostetal17b}.
The holomorphicity of \(T\) and \(J\) is somewhat surprising, but obtains a deep geometric explanation in the fully supersymmetric model below.

More generally, the scheme of converting a supersymmetric action principle of quantum field theory into a mathematical version with commuting fields works also in other cases, like the super Liouville action, and leads to a rich mathematical theory, see for instance \cite{Jost/Wang/Zhou07,Jost/Wang/Zhou/Zhu09,Jost/Wang/Zhou/Zhu14,Jost/Zhou/Zhu14,Jost/Zhou/Zhu15,Jost/Zhou/Zhu17a}.

This  seems to open  a beautiful research direction, with many difficult and challenging analytical problems, concerning the existence and regularity of solutions, and with the potential to profound geometric applications, possibly not unlike the Seiberg-Witten equations (see \cite{Seiberg/Witten94a,Seiberg/Witten94b} and for instance \cite{Jost17} and \cite{Jost96a} for the analytic aspects) or pseudoholomorphic curves giving rise to the Gromov-Witten invariants \cite{Gromov85,Taubes96b}.

\section{The supersymmetric action functional}%
\label{susy2}
We shall now switch from the Clifford algebra ${\CliffordAlgebraM}$ to ${\CliffordAlgebraP}$.
This will introduce a minus sign when we integrate by parts in \rf{intpar}.
In order to compensate, for this minus sign, we now choose $\psi$ as an \emph{anticommuting} field.
In physics terminology, $\psi$ then becomes a fermionic field, as it should.
And so will the gravitino field $\chi$. In particular, this will allow us to gain an additional symmetry, \emph{supersymmetry}. Mathematically, this means that we should enter the realm of supergeometry. This means that we can induce symmetries of the fields by transformations of the independent variables. Thus, we can capture supersymmetry as an extension or a partner of diffeomorphism invariance. And thus, the gravitino will become the superpartner of the metric, in the same manner as the spinor field $\psi$ will become a superpartner of the bosonic field $\phi$.

\subsection{The action functional}
In this section, we summarize results of \cite{Kessler17,Jost/Kessler/Tolksdorf16}.
The supersymmetric action functional that we want to investigate is
\begin{equation}
\label{AF}
	\begin{split}
		\ASRS(\phi, \psi;\gamma,\chi)&\coloneqq \int_\Sigma \left(\vphantom{\frac16}|\dd \phi|_{T^*\Sigma\otimes \phi^*TN}^2 + \langle \psi, \Dirac \psi \rangle_{\SpinorBundleP\otimes \phi^*TN} \right.\\
			&\hspace{4em} -4\langle (\mathrm{Id}\otimes\phi_*)(Q\chi), \psi \rangle_{\SpinorBundleP\otimes \phi^*TN} \\
			&\hspace{4em}\left.- |Q\chi|^2_{S\otimes TM} |\psi|^2_{\SpinorBundleP\otimes \phi^*TN}
			-\frac{1}{6} \Rm^{N}(\psi)\right) \dv_\gamma,
	\end{split}
\end{equation}
This functional was introduced in early works on string theory, see~\cite{DZ-CASS,BdVH-LSRIASS}, motivated by the search for a supersymmetric extension of the harmonic action functional.
In order to realize supersymmetry, physicists always work with anticommuting variables.
Indeed, it was argued in~\cite{KT-FSRSSCS} that, while the action~\(\ADHG\) can be motivated purely geometrically, the requirement of supersymmetry needs anticommuting spinors in general.

Consequently, while~\(\ASRS\) formally looks like \(\ADHG\) introduced in~\eqref{eq:AF}, we do now assume that~\(\psi\) and~\(\chi\) are anticommuting fields and we shall work with the Clifford-algebra \(\CliffordBundleP\) as in the physics literature.
Those two changes together assure that the Dirac action is nontrivial, compare Equation~\eqref{intpar}.

As \(\CliffordAlgebraP\) is isomorphic to the algebra of real \(2\times 2\)-matrices, the corresponding spinor bundle \(\SpinorBundleP=\SpinStructure\times\R^2\) is of real rank two.
The almost complex structure \(\ACRS\) of \(T\Sigma\) can be extended to \(\SpinorBundleP\) by setting
\begin{equation}
	\ACSBP s = -\omega s = -e_1\cdot e_2\cdot s
\end{equation}
for spinors \(s\).
The compatibility of the almost complex structures \(\ACSBP\) and \(\ACRS\) yields \(\ACSBP\left(X\cdot s\right) = \left(\ACRS X\right)\cdot s = - X\cdot \ACSBP s\) for all vectors \(X\) and spinors \(s\).
Consequently,
\begin{equation}
	\gamma\left(\Gamma(s,t), X\right) = \left< X\cdot s, t\right>_\SpinorBundleP
\end{equation}
defines a complex linear isomorphism of line bundles \(\Gamma\colon \SpinorBundleP\otimes\SpinorBundleP\to T\Sigma\).
From now%
\footnote{Here we use different convention for the gravitino than in~\cite{Jost/Kessler/Tolksdorf16}.
	There a gravitino is a section of \(T^*\Sigma\otimes S\).
	Furthermore the complex structure \(\ACSBP\) differs from the one in~\cite{Jostetal17} by a minus sign, but agrees with the one in~\cite{Jost/Kessler/Tolksdorf16}.
	It is this sign that leads to the isomorphism \(\SpinorBundleP\otimes_\C\SpinorBundleP=T\Sigma\) instead of \(\SpinorBundleP\otimes_\C\SpinorBundleP=T^*\Sigma\).
}%
on we assume that \(\psi\) is a section of \(\SpinorBundleP\otimes\phi^*TN\) and \(\chi\) a section of \(\SpinorBundleP\otimes T\Sigma\).
The bundle \(\SpinorBundleP\otimes T\Sigma\) decomposes into irreducible \(\SpinGroup\)-representations as \(\SpinorBundleP\oplus\SpinorBundleP\otimes_\C\SpinorBundleP\otimes_\C\SpinorBundleP\).
We denote again by \(Q\) the projection on the \(\frac32\)-part.

Despite the change of conventions for the Clifford-algebra the functional \(\ASRS\) shares the invariances of \(\ADHG\).
The fact that we work with anticommuting spinor fields yields an additional invariance that gives rise to the so-called supersymmetry of the action.
Thus, let us list the invariances of \rf{AF}.
\begin{enumerate}
	\item Conformal transformations: \(\ASRS(\phi,e^{-u}\psi;e^{2u}\gamma,e^{-2u}\chi) =\ASRS(\phi,\psi;\gamma,\chi)\).
	\item Super Weyl transformations: \(A(\phi, \psi, \gamma , \chi + \chi') = A(\phi, \psi, \gamma , \chi)\) for \(Q\chi'=0\).
	\item Diffeomorphisms of \(\Sigma\): \(A(\phi\circ f, f^*\psi, f^*\gamma, f^*\chi) = A(\phi, \psi, \gamma, \chi)\).
	\item Supersymmetry, that is first order invariance under:
		\begin{equation}\label{eq:SusyTrafo}
			\begin{aligned}
				\delta\phi &= \langle q, \psi \rangle_\SpinorBundleP &
				\delta\psi &= e_\alpha\cdot q\otimes \left(e_\alpha\phi - \langle\psi, \chi^\alpha\rangle\right)e_\alpha q\\
				\delta e_\alpha &= -2\langle e_\alpha\otimes \left(e_\beta\cdot q\right), \chi\rangle e_\beta &
				\delta\chi &= {\left(\nabla^\SpinorBundleP q\right)}^\sharp
			\end{aligned}
		\end{equation}
		where \(q\) is a spinor, \(e_\alpha\) is a $\gamma$-orthonormal frame and \(\nabla^\SpinorBundleP\) a  spin connection with torsion.
\end{enumerate}

Of course, all these symmetries induce Noether currents.
In order to understand the action functional \rf{AF} and its symmetries in a systematic manner, we shall introduce the formalism of supergeometry, and more precisely, the notion of a super Riemann surface.
In particular the local transformations~\eqref{eq:SusyTrafo}, which look somewhat ad-hoc will be geometrically interpreted as particular diffeomorphisms of the underlying supermanifold.
We shall then see that the same deep relationship that we have described in Sections \ref{teich}, \ref{vari} between the action functional~\rf{eq080106} and the geometry of the moduli space of Riemann surfaces extends to the present setting.
In fact, the Noether currents for the symmetries of $A$ will become cotangent vectors to the moduli space of super Riemann surfaces, and we can give a description of the latter in terms of metrics and gravitinos by dividing out symmetries.

\subsection{Super Riemann surfaces}
\subsubsection{Super geometry}
We recall some notions from supergeometry.
A \emph{supermanifold} is a locally ringed space $(\Vert M \Vert, \cO_M)$ locally isomorphic to $\R^{m|n}=(\R^m,C^\infty(\R^m,\R)\otimes_\R \Lambda_n)$ where~$\Lambda_n$ is the real Grassmann algebra of $n$ generators.
The elements of $\Lambda_n$ anticommute and are called the odd coordinate directions.
A supermanifold of dimension $(m|0)$ is simply an ordinary manifold of dimension $m$, also called an even supermanifold, whereas one of dimension $(0|n)$ is purely odd.
We shall only be interested here in the case $(m|n)=(2|2)$, that is, where we have two even and two odd dimensions.
We shall also use a complex notation and write $\C^{1|1}$ with coordinates $z=x+\ima y$ (even) and $\vartheta =\vartheta^1 + \ima\vartheta^2$ (odd).
A function $\Phi\colon \R^{2|2} \to \R$ is then given by
\bel{sug1}
\Phi^\sharp r= f_0(x) +\eta^1 \eta^2 f_{12}(x) \qquad  (r\in \R)
\qe
for even $x=(x^1,x^2)$ and odd $\eta=(\eta^1,\eta^2)$.
Here, we already see the principle, that we consider ordinary (smooth) functions of the even variables $x$ and expand in the odd variables $\eta$.
We observe that because the coordinates of $\eta$ anticommute, and therefore for instance $\eta^1 \eta^1 =0$, there can be no further terms in the expansion~\rf{sug1}.
If we also want to have odd functions $f$ in the expansion,  we need to extend the base and consider, for some purely odd parameter domain $B$,
\begin{equation*}
	\begin{split}
		\Phi&\colon  \R^{2|2} \times B \to \R \times B\\
		\Phi^\# r &= f_0(x) +\eta^\mu f_\mu(x)+ \eta^1 \eta^2 f_{12}(x)
	\end{split}
\end{equation*}
where the $f_\mu$ are odd and involve odd parameters of $B$.
Colloquially speaking, we simply pull some odd parameters out of the hat whenever we need them. Thus, we leave $B$ unspecified.

We mostly consider a complex supermanifold of complex dimension $(1|1)$.
Coordinate transformations then are of the form
\begin{equation}%
\label{sug5}
	\begin{split}
		(z,\vartheta)&\mapsto (z',\vartheta')\\
		z'&=f(z), \\
		\vartheta'&=\vartheta g(z).
	\end{split}
\end{equation}
More generally, when we consider a family over base $B$, we can also use odd functions of $z$ and therefore consider the more general coordinate transformations
\begin{align}
	z'&=f(z) +\vartheta\xi(z), &
	\vartheta'&=\zeta(z) + \vartheta g(z)
\end{align}
with odd \(\xi\), \(\zeta\).
As explained, this is only meaningful if they involve an odd parameter from $B$.

For all families of supermanifolds there exists a family $\Smooth{M}=(\Vert M\Vert, \cO_{\Smooth{M}})$ of even supermanifolds over $B$ and an embedding
\bel{sug3} i:\Smooth{M} \to M
\qe
which is the identity on the underlying topological space.
We point out that over families, the embedding \(i\) is not unique.

We also have the Berezin integral, which is defined for sections of the Berezin bundle \(\mathrm{Ber} T^*M\), the super generalization of the determinant line bundle.
The Berezin integral takes values in the functions on \(B\):
\begin{equation}
\label{sug4}
		\int_M \colon \mathrm{Ber} T^\ast M \to \cO_B\\
\end{equation}
One possible definition of the Berezin integral is to use the local formula
\begin{equation}
		\int_{\R^{m|n}} f(x,\eta) [\dx^1\dotsm\dx^m\dd\eta^1\dotsm\dd\eta^n] = \int_{\R^{m|0}} f_\mathrm{top}(x,\eta) \dx^1\dotsm\dx^m,
\end{equation}
and to globalize with the help of a partition of unity.
Here $f_\mathrm{top}$ is the coefficient of $\eta^n\dotsm\eta^1$, that is, of the highest order term occurring in the expansion of the function $f$ in terms of the $\eta^i$ (recall again that they are anticommuting).
Thus, each $d\eta^i$ in the integral simply cancels the coefficient $\eta^i$ in the expansion.

\subsubsection{Super Riemann surfaces}
A super Riemann surface is a $(1|1)$-dimensional complex supermanifold with an additional structure, a maximally non-integrable distribution $\cD$ of the tangent bundle $TM$ of rank $(0|1)$, i.e., an isomorphism
\begin{equation}
\label{sug7}
	\begin{split}
		\cD \otimes \cD & \to \faktor{TM}{\cD}\\
		X\otimes Y &\mapsto \faktor{[X,Y]}{\cD}.
	\end{split}
\end{equation}
This is the description of LeBrun and Rothstein~\cite{LBR-MSRS}.
According to them, we can find local coordinates $(z,\vartheta)$ so that $\cD$ is locally generated by
\begin{equation}
	D=\frac{\partial}{\partial \vartheta }+ \vartheta  \frac{\partial}{\partial z}.
\end{equation}
$D$ satisfies the key property
\begin{equation}
	D^2 = DD =\frac{\partial}{\partial z}
\end{equation}
which easily follows from the anticommutation rules for the odd variables.
Thus, $D$ is a square root of the ordinary derivative.
Since the latter can be considered as the generator of spatial translations, that is, of a family of diffeomorphisms, $D$ generates odd translations, or a family of superdiffeomorphisms.
A general superdiffeomorphism then combines those two types of translations.

Under the coordinate transformation \rf{sug5}, \(D\) is transformed to
\bel{sug9}
D = \vartheta f'(z) \frac{\partial}{\partial z'} + g(z) \frac{\partial}{\partial \vartheta'},
\qe
and for this in order to be equal to \(D'\), we need
\bel{sug10}
f'(z)=g(z)^2,
\qe
because then we get \(\cD'=g(z)\left(\vartheta' \frac{\partial}{\partial z'} + \frac{\partial}{\partial \vartheta'}\right)\).
This means that $\vartheta$ transforms as an (odd) section of the spinor bundle $\SpinorBundleP=T\Smooth{M}^{1/2}$.
This in turn means that $M$ corresponds to $(\Smooth{M},\SpinorBundleP)$ where $\Smooth{M}$ is a Riemann surface and $\SpinorBundleP$ is a spinor bundle.
There exist $2^{2p}$ spin structures on a Riemann surface of genus $p$, one for each element of $H^1(\Smooth{M},\Z_2)$, i.e., a choice of sign of the square root along 1-cycles.
While this observation completely classifies trivial families of super Riemann surfaces, it is important to study non-trivial familes, in particular those where the base has odd dimensions.

An alternative description of super Riemann surfaces is due to Giddings and Nelson~\cite{GN-GSRS}.
For them, a super Riemann surface is a $(2|2)$-dimensional real supermanifold with a reduction of the structure group to
\bel{sug11}
G=\left\{ \begin{pmatrix} B^2 & A \\ 0 & B \end{pmatrix}, A,B \in \C \right\}\subset \mathrm{Gl}_\C(1|1) \subset \mathrm{Gl}_\R(2|2)
\qe
with certain integrability conditions. A further reduction is achieved by
\begin{equation*}%
\label{sug12}
	\begin{split}
		U(1)&\to G\\
		U &\mapsto \begin{pmatrix} U^2 & 0 \\ 0 & U \end{pmatrix}.
	\end{split}
\end{equation*}
From that, we obtain a non-degenerate supersymmetric bilinear form on $TM$, given in $U(1)$-frames by
\begin{align}%
\label{sug13}
	m(F_a,F_b)&=\delta_{ab}, &
	m(F_a,F_\beta)&=0, &
	m(F_\alpha,F_\beta)&=\epsilon_{\alpha \beta},
\end{align}
where Latin indices $a,b$ stand for the even, Greek indices $\alpha, \beta$ for the odd directions. $\delta_{ab}$ is the symmetric Kronecker tensor, $\epsilon_{\alpha \beta} $ the corresponding antisymmetric one.

From \rf{sug7}, we get a split exact sequence
\begin{diag}%
\label{diag:SRSShortExactSequence}
	\matrix[seq](m) { 0 & \cD & TM=\cD^\perp\oplus \cD & \faktor{TM}{\cD} & 0.\\ };
	\path[pf]	(m-1-1) edge (m-1-2)
		(m-1-2) edge (m-1-3)
		(m-1-3) edge (m-1-4)
		(m-1-4) edge (m-1-5)
			edge[bend right=40] node[auto,swap]{\(p\)} (m-1-3);
\end{diag}
The differential of the embedding $i:\Smooth{M}\to M$ yields a different splitting of the pull back of \eqref{diag:SRSShortExactSequence} along \(i\), and so, we have
\begin{diag}
	\matrix[seq](m) { 0 & \SpinorBundleP & i^*TM & T\Smooth{M} & 0\\};
	\path[pf]	(m-1-1) edge (m-1-2)
		(m-1-2) edge (m-1-3)
		(m-1-3) edge (m-1-4)
		(m-1-4) edge (m-1-5)
			edge[bend right=40] node[auto,swap]{\(\tilde{p}\)} (m-1-3)
			edge[bend left=40] node[auto]{\(\dd{i}\)} (m-1-3);
\end{diag}
Here, $\tilde{p}$ induces $p_s:i^\ast TM \to \SpinorBundleP$.
We have the identifications
\bel{sug16}
\SpinorBundleP=i^\ast \cD, \qquad T\Smooth{M}=i^\ast \cD^\perp,
\qe
and the second equips $T\Smooth{M}$ with a metric $\gamma$, and the first makes $\SpinorBundleP$ a spinor bundle for that metric, because by \rf{sug7}
\begin{equation}%
\label{sug18}
	i^\ast \cD \otimes_\C i^\ast \cD =i^\ast \faktor{TM}{\cD} =T\Smooth{M}.
\end{equation}
The gravitino then is the section of $T\Smooth{M} \otimes \SpinorBundleP$ given by
\begin{equation}
\label{sug19}
	\gamma\left(\chi, v\right) = p_s(\tilde{p}-di)v.
\end{equation}
Thus, analogously to our description in Section \ref{teich} of a Riemannian metric as encoding the difference between the local Euclidean coordinate structure and the global conformal structure of a Riemann surface, we see here that the gravitino encodes the difference between the local coordinate structure realized by the embedding $i$ and the global super Riemann surface structure.
In particular, in the same manner that a metric can be locally made Euclidean by the uniformization theorem, a gravitino can be locally gauged to 0.
But, of course, not globally.
Of course, since the gravitino here is an odd field, regularity issues as for the uniformization theorem do not arise here.
On the other hand, from this perspective, we may ask what the appropriate regularity assumptions on the commuting gravitino field in Theorem \ref{theorem 1} might be. This depends on  whether also in that context, less regular fields can be transformed into more regular ones,
by exploiting the redundancy contained in the gravitino field. The point here is that since the Riemannian metric contains redundant information about the underlying conformal structure, we also expect that its superpartner, the gravitino, contains some redundant information.  We recall, for instance, that the action functional only involves the projection $Q\chi$ onto the spin $3/2$ component.

The metric, spinor bundle and gravitino entirely prescribe the super Riemann surfaces:
\begin{theo}\label{theo1}
	A $(2|0)$-dimensional Riemann surface $\Smooth{M}$ over a base $B$ with a spinor bundle $\SpinorBundleP$, a metric $\gamma$ and a gravitino $\chi$ determines a unique super Riemann surface $M$ with embedding $i:\Smooth{M}\to M$, such that the above construction gives $\gamma$ and $\chi$ back, up to Weyl (i.e., conformal) transformations of $\gamma$ and super Weyl transformations of $\chi$.
\end{theo}
The \emph{proof} depends on rather long and tedious computations employing a suitable adapted set of coordinates that were introduced by Wess and Zumino.

\begin{coro}
	There is a bijection
	\begin{equation}\label{sug20}
		\{ i:\Smooth{M}\to M, M \mathrm{ super\ Riemann\ surface}\} \longleftrightarrow \faktor{\{(\Smooth{M},\SpinorBundleP,\gamma,\chi)\}}{\mathrm{(super)Weyl}}
	\end{equation}
\end{coro}
This corollary is useful for deformation theory, because the right hand side involves no integrability conditions.

In fact, we even expect a global bijection
\begin{multline}%
\label{sug21}
	\faktor{\{ M \text{ super Riemann surface}\}}{\mathrm{(super) \;diffeomorphisms}} \quad\longleftrightarrow \\
	\faktor{\{(\Smooth{M},\SpinorBundleP,\gamma,\chi)\}}{{\mathrm{((super)\;Weyl, diffeomorphisms, supersymmetries)}}}
\end{multline}
where the superdiffeomorphisms induce both the ordinary diffeomorphisms of $\Smooth{M}$ and the supersymmetry transformations.
For a description of the super Teichmüller theory as a quotient of the superconformal structures up to superdiffeomorphisms, see Sachse~\cite{S-GAASTS}.

An infinitesimal variant of this expected global bijection is proven with the help of deformation theory:
\begin{theo}\label{theo2}
	The space of infinitesimal deformations of a super Riemann surface $M$ with embedding $i:\Smooth{M}\to M$ and $\gamma, \chi=0$ on $\Smooth{M}$ as above is
\bel{sug22}
H^0(T^\ast M \otimes_\C T^\ast M)\oplus H^0(\SpinorBundleP^\ast \otimes_C \SpinorBundleP^\ast \otimes_C \SpinorBundleP^\ast ).
\qe
\end{theo}
The proof decomposes variations \((h,r)\) of \((\gamma, \chi)\) as
\begin{align}
	\label{sug23}
	h &= e^u \gamma + L_X \gamma + \mathrm{susy}_q\gamma + \eta\\
	\label{sug24}
	r &= \chi' + L_X \chi + \mathrm{susy}_q\chi + \rho,
\end{align}
where \((e^u\gamma, \chi')\) are infinitesimal (super) Weyl transformations, $L_X$ is the Lie derivative in the direction of the vector field $X$, and \((\mathrm{susy}_q\gamma, \mathrm{susy}_q\chi)\) is the supersymmetry transformation of \((\gamma, \chi)\).
The remaining terms \((\eta, \rho)\) are the true deformations of the super Riemann surface.

\subsection{The symmetries of the action functional}
We have already listed the symmetries of the action functional~\eqref{AF} above. With the geometric concepts developed, we can now write the action functional in a manner that makes the both the structure of this functional and the analogies with the functional \rf{eq080106} transparent. In fact,
we can write
\begin{equation}%
\label{act2}
	\ASRS(\phi,\psi,\gamma,\chi) + \int_{\Smooth{M}} F^2 \dv_\gamma =
	\int_M \|\left.d\Phi\right|_\cD\|^2_{m\otimes \Phi^*g} [\dv_m],
\end{equation}
where the integral on the right is a Berezinian \rf{sug4}, and the field $\Phi$ contains both $\phi$ and $\psi$ as terms in its expansion with respect to the already mentioned Wess--Zumino coordinates
\begin{equation}
	\Phi=\phi +\eta^\mu \psi_\mu +\eta^1 \eta^2 F.
\end{equation}
The Euler--Lagrange equations for the field $F$ are  $F=0$, that is, $F$ becomes zero on-shell.
As the field \(F\) appears in the expansion of \(\Phi\), it is needed for the full superdiffeomorphism invariance of \(\ASRS\).
However, we do not enter this and simply omit the terms involving \(F\).

It should be noted that in~\eqref{act2} the action takes values in \(\cO_B\), hence a supercommutative ring possibly with nilpotent elements.
This is necessary to allow a mathematical treatment of the anticommuting fields \(\psi\) and \(\chi\).
As a consequence, the adaption of the analytical result for this variational problem requires further development of algebraic techniques.

The Berezin integral on the right-hand of~\eqref{act2} does not only look surprisingly similar to the harmonic action~\eqref{eq080106}, but also its supergeometric properties resemble those of the harmonic action functional.
Indeed, the Berezin integral does not depend on the chosen \(U(1)\)-reduction given by \(m\), but rather only on the super Riemann surface structure.
Every \(U(1)\)-reduction of the super Riemann surface will lead to the same value of the integral.
This is called superconformal invariance of the integral.
The superconformal invariance of the Berezin integral translates into the (super) Weyl invariance of \(\ASRS\).
Furthermore, the Berezin integral is superdiffeomorphism invariant.
Every superdiffeomorphism of~\(M\) can be decomposed into a diffeomorphism of \(\Smooth{M}\) and a diffeomorphism of~\(M\) that leaves \(i\colon\Smooth{M}\to M\) invariant.
The latter ones induce supersymmetry transformations on the fields \(\phi\), \(\psi\), \(\gamma\) and \(\chi\).
Hence the supersymmetry invariance of \(\ASRS\).

The superconformal and superdiffeomorphism invariance of~\(\ASRS\) indicate that~\(\ASRS\) should be a suitable tool to study the quotient of superconformal structures up to superdiffeomorphisms, i.e.\ the supermoduli space as in~\eqref{sug21}.
This would be in parallel to the harmonic maps approach to ordinary Teichmüller theory.

A first step in this direction is to study the Noether currents associated to the symmetries.
Recall the \emph{energy-momentum tensor} \(T\) and \emph{supercurrent} \(J\) defined by
\begin{align}
	\partial_\gamma \ASRS(\phi,\psi,\gamma,\chi)(\delta \gamma)&=\int_{\Smooth{M}} \left<\delta \gamma, T\right> \dv_\gamma, \\
	\partial_\chi \ASRS(\phi,\psi,\gamma,\chi)(\delta \chi)&= \int_{\Smooth{M}} \left<\delta\chi, J\right> \dv_\gamma.
\end{align}
As in Section~\ref{gravi1}, the Weyl invariance prescribes the trace of \(T\) and the super Weyl invariance prescribes \(J=QJ\).
Using the Euler--Lagrange equations for \(\phi\) and \(\psi\), the diffeomorphism invariance of \(\ASRS\) yields a coupled differential equation for \(T\) and \(J\) of divergence-type.
A second coupled equation for \(T\) and \(J\) is obtained from supersymmetry.

In the case \(\chi=0\) the equations decouple and give \(\diverg_\gamma T=0\) and \(\diverg_\gamma J=0\).
In that case, being trace and divergence free, $T$ can be considered as a holomorphic section of $(T^\ast M)^{\otimes 2}$, and the metric dual of $J$ as a holomorphic section of $\SpinorBundleP^*\otimes_\C \SpinorBundleP^*\otimes_\C\SpinorBundleP^*$.
Thus, according to Theorem \ref{theo2}, they constitute tangent vectors to the moduli space of super Riemann surfaces.

In conclusion, the action functional $A$ reflects the structure of the moduli space of super Riemann surfaces by its  symmetries and invariances in such a manner that its associated Noether currents become cotangent vectors to that space.
Thus, we can represent the (infinitesimal) geometry of that moduli space analytically through a physical action functional.
In turn, in order to understand the symmetries of that action functional, we can represent them as geometric variations of a superset of independent variables.
A further understanding of the critical points is expected to lead to deep insights on the Teichmüller theory of super Riemann surfaces.

\bibliographystyle{plain}
\bibliography{JKTWZ-HarmonicMapsToNLSMinQFT}

\end{document}